\documentclass[a4paper,12pt,leqno]{amsart}

\usepackage{latexsym}
\usepackage[all]{xy}

\usepackage{amsmath, amsfonts, amssymb,amscd,graphics,graphicx,color,eucal,setspace,xypic} 
\usepackage{comment}
\usepackage{mathtools}

\usepackage{graphicx}
\usepackage{stmaryrd}
\usepackage{bigdelim, multirow}
\usepackage[all]{xy}

\usepackage{ytableau} 

\usepackage{tikz}
\usepackage{tikz-cd}
\usepackage{tikz-3dplot}
\usetikzlibrary{decorations.pathmorphing}

\definecolor{gray}{gray}{0.7}
\definecolor{Gray}{gray}{0.3}

\usepackage{amscd}

\numberwithin{equation}{section}

\theoremstyle{break}
 \newtheorem{theorem}{Theorem}[section]
 \newtheorem{proposition}[theorem]{Proposition}
 
 \newtheorem{lemma}[theorem]{Lemma}

 \theoremstyle{definition}
 \newtheorem{definition}[theorem]{Definition}
 \newtheorem{remark}[theorem]{Remark}
 \newtheorem{example}[theorem]{Example}

\allowdisplaybreaks[4]

\def\C{\mathbb C}
\def\R{\mathbb R}

\def\Z{\mathbb Z}

\def\J{\mathcal{J}}
\def\L{\mathcal{L}}
\def\S{\mathfrak{S}}

\def\x{\mathbf{x}}

\def\hs{\widehat{\S_n}}
\def\ts{\widetilde{\S_n}}
\def\hwk{\widehat{W^K}}
\def\twk{\widetilde{W^K}}

\def\kostka{\mathcal{K}}

\def\rsk{\mathbf{R}\mathbf{S}\mathbf{K}}
\def\ssyt{\mathcal{S}\mathcal{S}\mathcal{Y}\mathcal{T}}
\def\syt{\mathcal{S}\mathcal{Y}\mathcal{T}}
\def\ch{\mathsf{c}\mathsf{h}}

\def\pr{\mathsf{p}\mathsf{a}\mathsf{r}}

\DeclareMathOperator{\pk}{pk}
\DeclareMathOperator{\Hop}{Hop}
\DeclareMathOperator{\des}{des}
\DeclareMathOperator{\Des}{Des}

\newcommand{\red}[1]{\textcolor{red}{#1}}

\begin{document}

\title[Gamma vectors of partitioned permutohedra]{Gamma vectors of partitioned permutohedra}

\author[T. Horiguchi]{Tatsuya Horiguchi}
\address{National Institute of Technology, Akashi College, 
Hyogo 674-8501, Japan}
\email{tatsuya.horiguchi0103@gmail.com}

\author[M. Masuda]{Mikiya Masuda}
\address{Osaka Central Advanced Mathematical Institute, Osaka Metropolitan University, 
Osaka 558-8585, Japan; \newline
\indent International Laboratory of Algebraic Topology and Its Applications, HSE University, Moscow 109028, Russia
}
\email{mikiyamsd@gmail.com}

\author[T. Sato]{Takashi Sato}
\address{Osaka Central Advanced Mathematical Institute, Osaka Metropolitan University, 
Osaka 558-8585, Japan}
\email{00tkshst00@gmail.com}

\author[J. Shareshian]{John Shareshian}
\address{Department of Mathematics, Washington University, St Louis, MO 63130, USA}
\email{jshareshian@wustl.edu}

\author[J. Song]{Jongbaek Song}
\address{Department of Mathematics Education, Pusan National University, Busan 46241, Republic of Korea}
\email{jongbaek.song@pusan.ac.kr}



\maketitle

\begin{abstract}
We determine that $\gamma$-vectors of partitioned permutohedra, thereby generalizing a result of Foata and Sch\"utzenberger.  Our result is closely related to a result of Athanasiadis on the representation of the symmetric group on the cohomology of the permutohedral variety.  We explain how to derive Athanasiadis' result from ours and vice versa. 
\end{abstract}

%

\section{Introduction}\label{intro}
We continue the study of partitioned permutohedra that was initiated by some of the present authors in \cite{HMSS}.  We examine the $\gamma$-vectors of these polytopes, obtaining a generalization of a result of Foata and Sch\"utzenberger found in \cite{FoataSchutzenberger} and recovering a representation-theoretic result of Athanasiadis found in \cite{Athanasiadis}.

The symmetric group $\S_n$ acts on Euclidean space $\R^n$ by permuting coordinates.  This action preserves each affine hyperplane $E\subset \R^n$ consisting of vectors with constant coordinate sum. Fix $p=(p_1,\ldots,p_n) \in E$ such that the entries $p_i$ are pairwise distinct.  The {\it permutohedron} $P_n$ is the convex hull of the $\S_n$-orbit of $p$.  (The combinatorial type of $P_n$ does not depend on the choice of $p$.)
For $k \in [n-1]$, we define the halfspace
$$
H(k)^\leq\coloneqq\{(x_1,\ldots,x_n )\in E \mid  x_k \leq x_{k+1}\}.
$$
The {\it partitioned permutohedron} associated to $K \subseteq [n-1]$ is
$$
P_n(K)\coloneqq P_n \cap \bigcap_{k \in K} H(k)^\leq.
$$

It was shown in \cite{HMSS} that the polytopes $P_n(K)$ are all simple of dimension $n-1$.  Therefore, if we write $f_i$ for the number of $i$-dimensional faces of $P_n(K)$ and set
$$
h_{P_n(K)}(t)\coloneqq\sum_{i=0}^{n-1}h_it^i=\sum_{i=0}^{n-1}f_i(t-1)^i,
$$
then $h_{P_n(K)}(t)$ is palindromic, that is, $h_i=h_{n-1-i}$ for $0 \leq i \leq n-1$. (This is the content of the {\it Dehn--Sommerville equations}, see for example \cite[Section 9.2]{Grun}.) It follows that there exist uniquely determined $\gamma_{n,j,K} \in \Z$ ($0 \leq j \leq \lfloor \frac{n-1}{2} \rfloor$) such that
$$
h_{P_n(K)}(t)=\sum_{j=0}^{\lfloor \frac{n-1}{2} \rfloor}\gamma_{n,j,K}t^j(1+t)^{n-1-2j}.
$$
We show here that each $\gamma_{n,j,K}$ is non-negative by exhibiting a set $\widehat{W^K}(j) \subseteq \S_n$ such that $\left|\widehat{W^K}(j)\right|=\gamma_{n,j,K}$.  

More precisely, let $W_K$ be the subgroup of $\S_n$ generated by the set of transpositions $\{(k,k+1):k \in K\}$.  Each coset of $W_K$ in $\S_n$ contains a unique shortest representative with respect to the usual length function (see Section \ref{prelim} for definitions of statistics on $\S_n$).  Write $W^K$ for the set of such shortest representatives.  Let $\widehat{\S_n}$ be the set of all $w \in \S_n$ with no double descent and no final descent and set
$$
\widehat{W^K}\coloneqq\widehat{\S_n} \cap W^K.
$$
Write $\des(w)$ for the number of descents of $w \in \S_n$ and for $0 \leq j \leq \lfloor \frac{n-1}{2} \rfloor$, set
$$
\widehat{W^K}(j)\coloneqq\{w \in \widehat{W^K}\mid \des(w)=j\}.
$$

\begin{theorem} \label{gammavector}
For each $j \in \lfloor \frac{n-1}{2} \rfloor$,
$$
\gamma_{n,j,K}=|\hwk(j)|.
$$
Equivalently,
$$
h_{P_n(K)}(t)=\sum_{w \in \hwk}t^{\des(w)}(1+t)^{n-1-2\des(w)}.
$$
\end{theorem}
In the main result of this article (Theorem \ref{theorem:main}), we also verify that 
\[
\gamma_{n,j,K} = |\{w \in \widetilde{W^K}\mid \des(w)=j\}|,
\] 
where $\widetilde{W^K}$ denotes the permutations in $W^K$ with no double descent and no \emph{initial} descent.  

\begin{remark}
It is known that $h_{P_n}(t)$ is the Eulerian polynomial,
$$
h_{P_n}(t)=A_n(t)\coloneqq\sum_{w \in \S_n}t^{\des(w)}.
$$
When $K=\emptyset$, we have $W_K=\{1\}$ and $W^K=\S_n$.  In this case, Theorem \ref{gammavector} specializes to
$$
A_n(t)=\sum_{W \in \widehat{\S_n}}t^{\des(w)}(1+t)^{n-1-2\des(w)},
$$
which was proved by Foata and Sch\"utzenberger in \cite{FoataSchutzenberger}.
\end{remark}

\begin{remark}
It was shown in \cite{HMSS} that each $P_n(K)$ is a simple flag polytope.  It was conjectured by Gal in \cite{Gal} that when the $h$-polynomial of a simple flag $d$-polytope is expanded in terms of the basis $t^j(1+t)^{d-2j}$, each coefficient is nonnegative.  Theorem \ref{gammavector} confirms that Gal's Conjecture holds for partitioned permutohedra.
\end{remark}  

If the chosen point $p$ has integer coordinates then $P_n(K)$ is a simple rational polytope and thus determines a toric variety $X_n(K)$ of complex dimension $n-1$.  By work of Danilov and of Jurkiewicz (see \cite[Theorem 10.8 and Remark 10.9]{Danilov}),
$$
h_{P_n(K)}(t)=\sum_{i=0}^{n-1}\dim H^{2i}(X_n(K);\C)t^i.
$$
By construction,  $P_n$ admits an $\S_n$ action.  With $X_n=X_n(\emptyset)$, this action determines an $\S_n$ representation on each $H^{2i}(X_n;\C)$. Such representations have been studied closely, for example by Procesi in \cite{Procesi}, by Stanley in \cite{Stanley1}, by Stembridge in \cite{Stembridge2,Stembridge1}, and by Dolgachev and Lunts in \cite{DolgachevLunts}. 

Theorem 1.1 of \cite{HMSS} says that the cohomology ring $H^\ast(X_n(K);\C)$ is isomorphic (as a graded ring) to the ring $H^\ast(X_n;\C)^{W_K}$ of $W_K$-invariant elements of $H^\ast(X_n;\C)$.  Every $\C[\S_n]$-module $V$ is completely determined by the fixed-point submodules $V^{W_K}$ as $K$ ranges over all subsets of $[n-1]$.  (Indeed, one needs only one representative from each $\S_n$-conjugacy class of parabolic subgroups.)  Since each $h_{P_n(K)}(t)$ is palindromic, we recover from Theorem \ref{gammavector} the already-known isomorphism of $\C[\S_n]$-modules
$$
H^{2i}(X_n;\C) \cong_{\S_n} H^{2(n-1-i)}(X_n;\C)
$$
for $0 \leq i \leq n-1$.  (See \cite[Corollary 4.4]{Stanley1} and  \cite[Proposition 12]{Stembridge1}.)  We should therefore hope for a nice expansion of the representation $\S_n$ on $H^\ast(X_n;\C)$ using the palindromic polynomials $t^j(1+t)^{n-1-2j}$.  Such an expansion was given by Athanasiadis, who expressed the irreducible decomposition of this representation in the desired form in \cite{Athanasiadis}.  (We remark here that Stembridge gave in \cite[Theorem 4.2]{Stembridge1} a combinatorial formula for the irreducible decomposition of each $H^{2k}(X_n;\C)$.  A different such formula is a special case of \cite[Theorem 6.3]{SWCQF}.)

Write $\syt_n$ for the set of all standard Young tableaux with $n$ boxes and arbitrary shape and $\widehat{\syt_n}$ for the set of all such tableaux having no double descent and no descent at $n-1$.  Given a tableau $Q$, write $\lambda(Q)$ for the shape of $Q$.  (See Section \ref{Athanasiadis' result} for this and other definitions.)  With $R_{n,j}$ denoting the $\C[\S_n]$-module $H^{2j}(X_n;\C)$, Athanasiadis proved the following theorem. 

\begin{theorem}(\mbox{Athanasiadis, see }\cite[Corollary 2.41, Equation (90)]{Athanasiadis}) \label{gammarep}
For $0 \leq j \leq \lfloor \frac{n-1}{2} \rfloor$,
$$
R_{n,j} \cong_{\S_n} \sum_{\substack{Q \in \widehat{\syt_n} \\ {\des(Q)=j}}}S^{\lambda(Q)}.
$$
Equivalently,
$$
\sum_{k=0}^{n-1}H^{2k}(X_n;\C)t^k\cong_{\S_n} \sum_{Q \in \widehat{\syt_n}}S^{\lambda(Q)}t^{\des(Q)}(1+t)^{n-1-2\des(Q)}.
$$
\end{theorem}

As we have discussed, Theorems \ref{gammavector} and \ref{gammarep} carry the same information; if one knows the dimensions of the fixed point spaces for all parabolic subgroups, then one knows the representation and the converse is clear.  After proving Theorem \ref{gammavector} directly (without using Theorem \ref{gammarep}) in Section \ref{sec_proof_of_main_thm}, we show in Section \ref{relation} how to derive each of the theorems quickly assuming the other.

\section*{Acknowledgements}
The authors are grateful to Christos A. Athanasiadis for his helpful comments on an earlier version of this paper. This research was partly supported by Osaka Central Advanced Mathematical Institute (MEXT Joint Usage/Research Center on Mathematics and Theoretical Physics). The first author is supported in part by JSPS Grant-in-Aid for Young Scientists: 19K14508. The second author was supported in part by Grant-in-Aid for Scientific Research 22K03292 and the HSE University Basic Research Program.  The fourth author was supported in part by National Science Foundation grant DMS-1518389. The last named author was supported by the National Research Foundation of Korea(NRF) grant funded by the Korea government(MSIT) (RS-2025-00555914)

\section{Preliminaries} \label{prelim}

We review here some standard definitions and known facts from algebraic combinatorics.  We also introduce some notation. It is reasonable to skip much of this section and refer back as necessary, but note the use of $\widehat{X}$ (resp. $\widetilde{X}$) to indicate the set of all elements of $X$ having no double descent and no final descent (resp. no initial descent) when $X$ is a set on whose elements descents are defined.

\subsection{Words, permutations, and descents}
By a {\it word of length $n$} we mean an ordered list $v=v(1) \cdots v(n)$ of $n$ positive integers.  So, the set of words of length $n$ is $\Z_{>0}^n$.  A {\it descent} of $v \in \Z_{>0}^n$ is any $i \in [n-1]$ such that $v(i)>v(i+1)$.  We write $\Des(v)$ for the set of descents of $v$ and $\des(v)$ for $|\Des(v)|$.  We say that $v$ has a {\it double descent} if there is some $i \in [n-2]$ such that $\{i,i+1\} \subseteq \Des(v)$ and that $v$ has an {\it initial descent} if $1 \in \Des(v)$.
Similarly, we say that $v$ has a {\it final descent} if $n-1 \in \Des(v)$. The {\it multiplicity} of $i \in \Z_{>0}$ in $v$ is the number of $j \in [n]$ such that $v(j)=i$, which will be denoted by $m_i(v)$.  For example, $v=4223185$ is a word of length seven, $\Des(v)=\{1,4,6\}$, $\des(v)=3$, and (for example) $m_2(v)=2$ and $m_8(v)=1$.  In this case, $v$ has no double descent but does have an initial descent and a final descent.  We write $\widetilde{\Z_{>0}^n}$ (resp. $\widehat{\Z_{>0}^n}$) for the set of all words of length $n$ having no double descent and no initial descent (resp. final descent).

For a subset $X$ of $\Z_{>0}^n$, we write $\widetilde{X}$ and $\widehat{X}$ for $X \cap \widetilde{\Z_{>0}^n}$ and $X \cap \widehat{\Z_{>0}^n}$, respectively. 
For example, a permutation $w \in \S_n$, written in one-line notation, can be thought of as a word of length $n$ in which $m_i(v)=1$ for all $i \in [n]$.  
Then $\ts$ is the set of permutations with no double descent and no initial descent.

\subsection{Parabolic subgroups of $\S_n$}
A subgroup of $\S_n$ is {\it parabolic} if it is generated by transpositions.  Every parabolic subgroup is conjugate in $\S_n$ to a {\it standard parabolic subgroup}, which is a subgroup generated by some adjacent transpositions $(k,k+1)$.  For a subset $K \subseteq [n-1]$, we write $W_K$ for the standard parabolic subgroup generated by $\{(k,k+1):k \in K\}$.  If $1 \leq i<j \leq n$ then $i,j$ lie in the same $W_K$-orbit on $[n]$ if and only if $\ell \in K$ whenever $i \leq \ell<j$. 

The {\it length} $\ell(w)$ of $w \in \S_n$ is the smallest number of adjacent transpositions whose product is $w$, which is equal to the number of pairs $(i,j)$ from $[n]$ such that $i<j$ and $w(i)>w(j)$ (see for example \cite[Proposition 1.5.2]{BB}).  It is known (see for example \cite[Proposition 2.4.4]{BB}) that for $K \subseteq [n-1]$, each left coset $C$ of $W_K$ in $\S_n$ contains some $w$ such that $\ell(w)<\ell(x)$ for all $x \in C \setminus \{w\}$. We compose permutations from left to right.
So, applying $x \in wW_K$ to $[n]$ is the same as applying $w$ and then permuting elements within $W_K$-orbits.  As the length of $w$ is minimal in $wW_K$, it must be the case that $i$ appears to the left of $j$ in $w$ whenever $i<j$ and $i,j$ are in the same $W_K$-orbit.  We write $W^K$ for the set of all minimal length coset representatives for $W_K$ in $\S_n$.
In other words, if we decompose a subset $K \subseteq [n-1]$ into maximal consecutive strings $K=[a_1,b_1] \sqcup \cdots \sqcup [a_m,b_m]$, then $W^K$ is the set of elements of $\mathfrak{S}_n$ satisfying the following condition: 
\begin{align}
\begin{split} \label{eq:KW}
&\textrm{for each $i$ with $1 \leq i \leq m$, the numbers $a_i, a_i+1, \cdots,b_i,  b_i+1$}  \\
&\textrm{appear in the one-line notation of $w$ as a subsequence in increasing order.} 
\end{split}
\end{align}

\subsection{Fans, polytopes, toric varieties, and $h$-polynomials}
We will not discuss at length the rich theory of toric varieties and related combinatorial objects, which can be found in various textbooks, including \cite{Oda}.  Rather, we record a little more than we need herein.  A $d$-dimensional convex polytope $P$ in $\R^d$ is {\it simple} if each vertex of $P$ is contained in exactly $d$ faces of dimension $(d-1)$.  Given a simple  polytope $P$, we write $f_i=f_i(P)$ for the number of $i$-dimensional faces of $P$ ($0 \leq i \leq d$).  The {\it $f$-polynomial} of $P$ is
$$
f_P(t)\coloneqq\sum_{i=0}^d f_it^i,
$$ 
and the $h$-polynomial of $P$ is
$$
h_P(t)\coloneqq f_P(t-1).
$$
So, for example, if $P$ is a $4$-gon in $\R^2$ then $f_P(t)=4+4t+t^2$ and $h_P(t)=1+2t+t^2$.

Given a simple polytope $P \subseteq \R^n$ with rational vertices, the fan $\Delta_P$ dual to $P$ is simplicial and determines a toric variety $X_P$.  The key point here, mentioned in the introduction, is the result proved in the projective case by Jurkiewicz and in the general case by Danilov,
$$
\sum_{k=0}^d \dim H^{2k}(X_P;\C)t^k = h_P(t).
$$
We remark that $X_P$ has trivial rational cohomology in odd degrees, hence $h_P(t)$ tells us everything about the Betti numbers of $X_P$.

\subsection{Partitioned permutohedra}
Let $p_1,\ldots,p_n$ be pairwise distinct real numbers.  The {\it permutohedron} $P_n(p)$ associated to $p=(p_1,\ldots,p_n)$ is the convex hull of the $n!$ points $(p_{w(1)},\ldots,p_{w(n)})$, with $w$ running through $\S_n$. We note that $P_n(p)$ lies in the affine hyperplane $E \subset \R^n$ consisting of vectors whose entries sum to $\sum_{i=1}^n p_i$. It is also known that $P_n(p)$ is an $(n-1)$-dimensional simple polytope.  Under our assumption that the coordinates in $p$ are pairwise distinct, the combinatorial type of $P_n(p)$ does not depend on the choice of $p$.  Therefore, we write $P_n$ for $P_n(p)$ for simplicity.  The toric variety associated with $P_n$ is called the {\it permutohedral variety}.

For each $k \in [n-1]$, we set
\[
H(k)^\leq \coloneqq \{(x_1,\ldots,x_n) \in E\mid  x _k \leq x_{k+1}\}.
\]
The {\it partitioned permutohedron} $P_n(K)$ associated to $K \subseteq [n-1]$ is
\[
P_n(K):=P_n \cap \bigcap_{k \in K}H(k)^{\leq}.
\]
The polytopes $P_n(K)$ were first studied in \cite{HMSS}.  It is known that, in addition to being simple, $P_n(K)$ is a flag polytope (see \cite[Proposition 4.6]{HMSS}).  Note that $P_n(\emptyset)=P_n$.  On the other hand, $P_n([n-1])$ is combinatorially equivalent to a cube.  Using work of Precup on Hessenberg varieties in \cite{Precup}, the following result was proved in \cite{HMSS}.

\begin{proposition}[\cite{HMSS}, Proposition 7.4] 
\label{proposition:h-polynomial PnK}
Given $K \subseteq [n-1]$, set
$$
W(K):=\{w \in \S_n\mid w^{-1}(k)-w^{-1}(k+1) \leq 1 \mbox{ for all } k \in K\}.
$$
For each such $K$,
$$
h_{P_n(K)}(t)=\sum_{w \in W(K)}t^{\des(w)}.
$$
\end{proposition}

\section{Main theorem}\label{sec_main}
The main result of this paper gives a combinatorial formula for expanding $h_{P_n(K)}(t)$ in terms of the gamma basis $\{t^j (1+t)^{n-1-2j} \}_{0\leq j \leq \lfloor \frac{n-1}{2} \rfloor}$. 
Foata and Sch\"{u}tzenberger  proved the following result. 

\begin{theorem}[\cite{FoataSchutzenberger}] \label{theorem:gamma_vector_Pn}
We can write the Eulerian polynomial as 
\begin{align*}
A_n(t) = h_{P_n}(t) =  \sum_{j=0}^{\lfloor \frac{n-1}{2} \rfloor} \gamma_{n,j} \, t^j (1+t)^{n-1-2j},
\end{align*}
where the coefficient $\gamma_{n,j}$ is equal to
$$
\gamma_{n,j} = |\{w \in \ts \mid \des(w)=j \}| = |\{w \in \hs \mid \des(w)=j \}|.
$$
\end{theorem}
Note that the last equality above follows from the involution $\ts \to \hs$ sending $w$ to the permutation $w^\ast$ satisfying $w^\ast(i)=n+1-w(n+1-i)$ for all $i \in [n]$.
Our main theorem generalizes Theorem~\ref{theorem:gamma_vector_Pn} as follows.

\begin{theorem} \label{theorem:main}
For a subset $K \subseteq [n-1]$, the $h$-polynomial of the partitioned permutohedron $P_n(K)$ satisfies 
\begin{align*}
h_{P_n(K)}(t) =  \sum_{j=0}^{\lfloor \frac{n-1}{2} \rfloor} \gamma_{n,j,K} \, t^j (1+t)^{n-1-2j},
\end{align*}
where for each $j$,  
$$
\gamma_{n,j,K} = |\{w \in \twk \mid \des(w)=j \}| = |\{w \in \hwk \mid \des(w)=j \}|.
$$ 
\end{theorem}

\begin{example}
Consider the case $n=5$ and $K=\{1,3\}$.  The orbits of $W_K$ are $\{1,2\}$, $\{3,4\}$ and $\{5\}$.  So, $W^K$ consists of all $w \in \S_5$ such that $w^{-1}(1) \leq w^{-1}(2)$ and $w^{-1}(3) \leq w^{-1}(4)$.  Among the $30$ such $w$, there are $6$ satisfying $w(1)=5$ and thus having an initial descent.  There are another $6$ satisfying $w(1)=3$ and $w(2)=1$.  Among the remaining $18$ permutations, $13542$, $15324$, and $35412$ each have a double descent.  We examine the $15$ elements of $\widetilde{W^K}$ and conclude that
\begin{itemize}
\item $12345$ has no descent;
\item $12354$, $12534$, $13245$, $13425$, $13452$, $13524$, $15234$, $34125$, $34512$, and $35124$ all have one descent; and
\item $13254$, $15342$, $34152$, and $35142$ all have two descents.
\end{itemize}
Using Theorem \ref{theorem:main}, we conclude that
\begin{equation} \label{example3.3eq1}
h_{P_5(\{1,3\})}=(1+t)^4+10t(1+t)^2+4t^2=1+14t+30t^2+14t^3+t^4.
\end{equation}

We can confirm that \eqref{example3.3eq1} is correct using Proposition \ref{proposition:h-polynomial PnK}.  We observe that $W(K)$ consists of all $w \in \S_5$ such that when $w$ is written in one-line notation, both
\begin{itemize}
\item[(a)] either $2$ appears to the right of $1$ or $2$ appears immediately to the left of $1$, and
\item[(b)] either $4$ appears to the right of $3$ or $4$ appears immediately to the left of $3$
\end{itemize}
hold.  We list in the table below all fourteen~$v \in \S_4$ satisfying (a) and (b).  For each such~$v$, we let $\S_5(v)$ be the set of all $w \in \S_5$ that are obtained by inserting $5$ into the one-line representation of $v$ so that (a) and (b) continue to hold, and give $\sum_{w \in \S_5(v)}t^{\des(w)}$.
$$
\begin{array}{|c|c|}
\hline v & \sum_{w \in \S_5(v)}t^{\des(w)} \\ \hline 
1234 & 1+4t \\ \hline 
1324,1342,3124,3412 & 2t+3t^2 \\ \hline 
3142 & 3t^2+2t^3 \\ \hline 
3214,3421,1432,4312 & 2t^2+2t^3 \\ \hline 
2134,1243 & t+3t^2 \\ \hline 
2143 & t^2+2t^3 \\ \hline 
4321 & 2t^3+t^4 \\ \hline
\end{array}
$$
Direct calculation shows that \eqref{example3.3eq1} holds.
\end{example}

\section{Valley hopping}
In this section, we review the \emph{valley-hopping} operation on a permutation, defined by Foata and Strehl in \cite{FoataStrehl}, where they use this to give a proof of Theorem~\ref{theorem:gamma_vector_Pn}. See also the description by Petersen in \cite{Petersen}.

\begin{definition}
For a permutation $w=w(1)w(2)\cdots w(n) \in \mathfrak{S}_n$ in one-line notation, we say 
\begin{enumerate}
\item $w(i)$ is a \textit{peak} if $w(i-1)< w(i) > w(i+1)$.
\item $w(i)$ is a \textit{valley} if $w(i-1)> w(i) < w(i+1)$.
\item $w(i)$ is a \textit{free letter} if $w(i)$ is neither peak nor valley.
\end{enumerate}
Here we take $w(0)=w(n+1)=\infty$  by convention, so that neither $w(1)$ nor $w(n)$ is a peak. 
\end{definition}

\begin{example}
If $w=672841359 \in \mathfrak{S}_9$, then $7,8$ are peaks, $6,2,1$ are valleys, and $4,3,5,9$ are free letters. We draw a permutation as a \emph{mountain range} following \cite[Section 4]{Petersen}. 
See Figure~\ref{picture:peaks, valleys, and free letters.} for the description of $w$. 
\begin{figure}
\begin{tikzpicture}[scale=0.5]
\draw (0,10) node[inner sep=1pt, draw, fill=white, circle] {\footnotesize$\infty$}
--(1,6) node[inner sep=1pt, draw,fill=white, circle] {\footnotesize$6$}
--(2,7) node[inner sep=1pt, draw,fill=white, circle] {\footnotesize$7$}
-- (3,2) node[inner sep=1pt, draw,fill=white, circle] {\footnotesize$2$}
-- (4,8) node[inner sep=1pt, draw,fill=white, circle] {\footnotesize$8$}
-- (5,4) node[inner sep=1pt, draw,fill=white, circle] {\footnotesize$4$}
-- (6,1) node[inner sep=1pt, draw,fill=white, circle] {\footnotesize$1$}
-- (7,3) node[inner sep=1pt, draw,fill=white, circle] {\footnotesize$3$}
-- (8,5) node[inner sep=1pt, draw,fill=white, circle] {\footnotesize$5$}
-- (9,9) node[inner sep=1pt, draw,fill=white, circle] {\footnotesize$9$}
-- (10,10) node[inner sep=1pt, draw,fill=white, circle] {\footnotesize$\infty$};
\end{tikzpicture}
\caption{Peaks, valleys and free letters for $w=672841359$.}
\label{picture:peaks, valleys, and free letters.}
\end{figure}
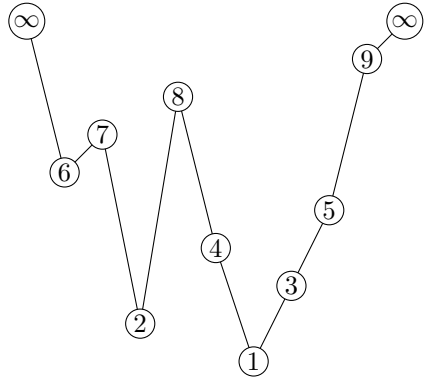
\end{example}

Let $w(i)=j$ be a free letter. We denote by $H_j(w)$  the permutation obtained by moving~$j$ directly across the adjacent valley(s) to the nearest mountain slope of the same height. More precisely, we define $H_j(w)$ as follows.
\begin{enumerate}
\item[$\bullet$] If $w(i) = j$ lies on a downslope, i.e., $w(i-1) > w(i) > w(i+1)$, then we define 
$$
H_j(w) = w(1) \cdots w(i-1)w(i+1) \cdots w(k) \, j \, w(k+1) \cdots w(n)
$$
for the smallest $k > i$ such that $w(k) < j < w(k+1)$.
\item[$\bullet$] If $w(i) = j$ lies on an upslope, i.e., $w(i-1) < w(i) < w(i+1)$, then we define 
$$
H_j(w) = w(1) \cdots w(k-1) \, j \, w(k) \cdots w(i-1)w(i+1) \cdots w(n)
$$
for the largest $k < i$ such that $w(k-1) >  j > w(k)$.
\end{enumerate}

\begin{example}
For $w=672841359\in \mathfrak{S}_9$, we have $H_4(w)=672813459$ and $H_9(w)=967284135$. 
See Figure \ref{fig_valley_hopping}. 
\begin{figure}
\begin{tikzpicture}[scale=0.5]
\draw[-stealth, red] (5.5,4)--(7.3,4);
\node[above] at (6.3,4) {\red{\footnotesize{move}}};
\draw (0,10) node[inner sep=1pt, draw, fill=white, circle] {\footnotesize$\infty$}
--(1,6) node[inner sep=1pt, draw,fill=white, circle] {\footnotesize$6$}
--(2,7) node[inner sep=1pt, draw,fill=white, circle] {\footnotesize$7$}
-- (3,2) node[inner sep=1pt, draw,fill=white, circle] {\footnotesize$2$}
-- (4,8) node[inner sep=1pt, draw,fill=white, circle] {\footnotesize$8$}
-- (5,4) node[inner sep=1pt, draw, fill=white, circle] {\footnotesize$4$}
-- (6,1) node[inner sep=1pt, draw,fill=white, circle] {\footnotesize$1$}
-- (7,3) node[inner sep=1pt, draw,fill=white, circle] {\footnotesize$3$}
-- (8,5) node[inner sep=1pt, draw,fill=white, circle] {\footnotesize$5$}
-- (9,9) node[inner sep=1pt, draw,fill=white, circle] {\footnotesize$9$}
-- (10,10) node[inner sep=1pt, draw, fill=white, circle] {\footnotesize$\infty$};
\node at (5,-1) {$H_4(w)=672813459$};

\begin{scope}[xshift=400]
\draw[-stealth, red] (8.5,9)--(0.4,9);
\node[above] at (5,9) {\red{\footnotesize{move}}};
\draw (0,10) node[inner sep=1pt, draw, fill=white, circle] {\footnotesize$\infty$}
--(1,6) node[inner sep=1pt, draw,fill=white, circle] {\footnotesize$6$}
--(2,7) node[inner sep=1pt, draw,fill=white, circle] {\footnotesize$7$}
-- (3,2) node[inner sep=1pt, draw,fill=white, circle] {\footnotesize$2$}
-- (4,8) node[inner sep=1pt, draw,fill=white, circle] {\footnotesize$8$}
-- (5,4) node[inner sep=1pt, draw, fill=white, circle] {\footnotesize$4$}
-- (6,1) node[inner sep=1pt, draw,fill=white, circle] {\footnotesize$1$}
-- (7,3) node[inner sep=1pt, draw,fill=white, circle] {\footnotesize$3$}
-- (8,5) node[inner sep=1pt, draw,fill=white, circle] {\footnotesize$5$}
-- (9,9) node[inner sep=1pt, draw ,fill=white, circle] {\footnotesize$9$}
-- (10,10) node[inner sep=1pt, draw,fill=white, circle] {\footnotesize$\infty$};
\node at (5,-1) {$H_9(w)=967284135.$};
\end{scope}
\end{tikzpicture}
\caption{Valley-hopping.}
\label{fig_valley_hopping}
\end{figure}
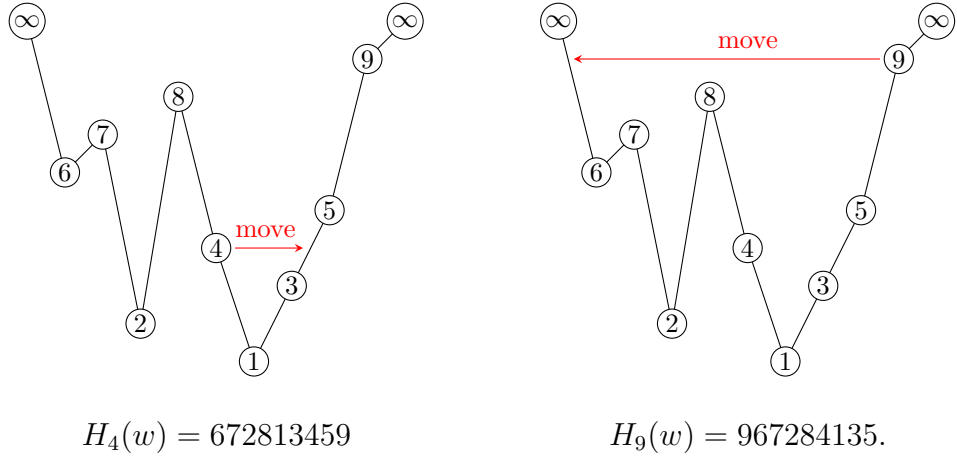
\end{example}

For any collection of free letters $J = \{j_1, \ldots , j_k\}$, we define 
$$
H_J(w) \coloneqq H_{j_1} \cdots H_{j_k}(w)
$$
and the hop equivalence relation $\sim$ on $\mathfrak{S}_n$ by
$$
w \sim u \text{ if and only if } u = H_J(w) 
\text{ for some collection of free letters } J = \{j_1, \ldots , j_k\}.
$$ 
We denote by $\Hop(w)$ the hop equivalence class of $w$. We observe that $H_J(w)$ has the same peaks and valleys as $w$.

%
%
%
%
%

For any $w \in \mathfrak{S}_n$, one can see that 
\begin{align} \label{eq:descent}
\des(w) = |\{ \textrm{peaks of } w \}| + |\{ \textrm{free letters of } w \textrm{ which lie on a downslope} \}|.
\end{align}
If $w$ has $r$ peaks, then it has $r+1$ valleys and $n-1-2r$ free letters. So, writing $\pk(w)$ as the number of peaks of $w\in \mathfrak{S}_n$, 
we have 
\begin{align} \label{eq:2-1}
\sum_{u \in \Hop(w)} t^{\des(u)} = t^{\pk(w)} (1+t)^{n-1-2\pk(w)}. 
\end{align}
Indeed, if $u_0$ is a \emph{canonical representative} of $\Hop(w)$ obtained by putting free letters on upslopes, then we have $t^{\des(u_0)}=t^{\pk(u_0)}$. 
If $u \in \Hop(w)$ is obtained from $u_0$ by moving~$m$ free letters, then we have $t^{\des(u)}=t^{\pk(u)+m}$. 
Hence, the left-hand side of \eqref{eq:2-1} can be expressed as 
\begin{align*}
t^{\pk(u_0)} \cdot \underbrace{(1+t)(1+t)\cdots(1+t)}_{\# \textrm{of free letters}},
\end{align*}
as desired.

One can easily see that $\ts$,  the set of permutations with no double descent and no initial descent,  coincides with the set of the canonical representative above for each hop equivalence class. 
In other words,
\begin{align}\label{eq:hatSn} 
\begin{split}
\ts = &\{w \in \mathfrak{S}_n \mid \textrm{all free letters of $w$ lie on an upslope} \} \\
= &\{w \in \mathfrak{S}_n \mid \pk(w)=\des(w) \}  
\end{split}
\end{align}
where the second equality follows from \eqref{eq:descent}. 
By summing \eqref{eq:2-1} over all $w \in \ts$, we get
\begin{align*} 
\sum_{w \in \ts} \sum_{u \in \Hop(w)} t^{\des(u)} = \sum_{w \in \ts} t^{\pk(w)} (1+t)^{n-1-2\pk(w)}. 
\end{align*}
This gives a proof of Theorem~\ref{theorem:gamma_vector_Pn}. We refer to \cite[Chaper 4]{Petersen} for more details. 

Recall that for a subset $K \subseteq [n-1]$, we define 
\begin{align} \label{eq:hatKSn}
\twk =\{w \in \ts \mid w \ \textrm{satisfies} \ \eqref{eq:KW} \}.
\end{align}
By summing \eqref{eq:2-1} over all $w \in \twk$, we get
\begin{align} \label{eq:3-1} 
\sum_{w \in \twk} \sum_{u \in \Hop(w)} t^{\des(u)} = \sum_{w \in \twk} t^{\pk(w)} (1+t)^{n-1-2\pk(w)}. 
\end{align}
Hence, in order to prove Theorem~\ref{theorem:main},  it suffices to show that the left-hand side of \eqref{eq:3-1} is equal to $h_{P_n(K)}(t)$.

Recall from Proposition~\ref{proposition:h-polynomial PnK} that a combinatorial formula for $h_{P_n(K)}(t)$ is given by 
\begin{align} \label{eq:h-polynomial PnK}
h_{P_n(K)}(t) = \sum_{v \in W(K)} t^{\des(v)}.
\end{align}
Here, we note that $v\in  W(K)$ if and only 
if in the one line notation for $v$, we see either
\begin{equation} \label{eq_cond_WK}
\text{(i) } \cdots k \cdots k+1 \cdots \quad \text{ or }  \quad  \text{(ii) }\cdots k+1~k \cdots
\end{equation}
for each $k\in K$. 
Namely, the cases (i) and (ii) above mean $v^{-1}(k)<v^{-1}(k+1)$ and $v^{-1}(k) -1 = v^{-1}(k+1)$, respectively.

When $K=\emptyset$, namely for the case where $P_n(K)$ is a permutohedron, we have $W(K)=\mathfrak{S}_n$ and it can be decomposed into $\bigsqcup_{w \in \ts} \Hop(w)$, which does not hold for arbitrary $K\subseteq [n-1]$ in general. Instead, we claim that there is a bijection 
\begin{equation}\label{eq_bijection}
\Theta \colon \{(w,u) \in \twk \times \mathfrak{S}_n \mid u \in \Hop(w) \} \to W(K)
\end{equation}
defined by $(w,u)\mapsto v$ such that $\des(u)=\des(v)$. Then, the bijection $\Theta$ together with~\eqref{eq:h-polynomial PnK} shows that $h_{P_n(K)}(t)$ agrees with the left-hand side of \eqref{eq:3-1}, which is a key step to prove Theorem~\ref{theorem:main}. 
For the purpose of a constructing a bijection $\Theta$, we introduce the following operator.
\begin{definition} \label{definition:key operator}
\begin{enumerate}
\item For a permutation $u \in \mathfrak{S}_n$ and $k \in [n]$, if the one line notation of $u$ does not satisfy the condition \eqref{eq_cond_WK}, then we define $\J_k(u)$ to be the permutation obtained from $u$ by moving $k+1$ to the right so that it is of the form (ii) of \eqref{eq_cond_WK}. Otherwise, we define $\J_k(u)=u$.
\item In general, for $K=\{i_1, i_2, \dots, i_m\}$ with $i_1<i_2<\cdots<i_m$, we define 
\begin{align*}
\J_K(u) \coloneqq \J_{i_m} \cdots \J_{i_2} \J_{i_1} (u).
\end{align*}
We observe that $\J_K(u) \in W(K)$ for all $u \in \mathfrak{S}_n$ by the definition.
\end{enumerate}
\end{definition}

We will prove the following proposition in the next two sections.

\begin{proposition} \label{proposition:key}
For a subset $K \subseteq [n-1]$, the map $\Theta$ of \eqref{eq_bijection} 
defined by 
\[
\Theta(w,u) = \J_K(u)
\] 
satisfies $\des(\J_K(u))=\des(u)$ and is bijective. 
\end{proposition}

\begin{example}
When $n=3$ and $K=\{1, 2 \}$,  we have $\twk = \{123 \}$ and 
\[
\Hop(123)=\{123, 213, 312, 321\}.
\]
Applying the operator $\J_K$ to permutations in $\Hop(123)$, we obtain $\{123, 213, 132, 321\}$ which  agrees with $W(K)$.
\end{example}

\section{Preserving descent} \label{sec_Dec_Prev}
In this section we prove that the operator $\J_K$ introduced in Definition \ref{definition:key operator} preserves the number of descents, i.e., $\des (\J_K(u))=\des (u)$. 
We first prove the claim for the case where~$K$ consists of consecutive integers. 

\begin{lemma} \label{lemma:descent-preserving for closed interval}
Let $K=[k,k+\ell-1]$ be the closed interval from $k$ to $k+\ell-1$. 
Let $w \in \twk$ and $u \in \Hop(w)$. 
Then we have
\begin{align*}
\des(\J_K(u))=\des(u).
\end{align*}
\end{lemma}

\begin{proof}
We prove this by induction on the cardinality $|K|=\ell$.
The base case is $K = \{k \}$. 
Take an element $w \in \widetilde{W^{\{k\}}}$. 
Note that $k+1$ is northeast of $k$ in the mountain range description of $w$ by \eqref{eq:hatKSn}. 
Consider the situation $\J_k(u) \neq u$ for $u \in \Hop(w)$, which happens only in the following two cases.

\smallskip 
\noindent
\textbf{Case (a):} The numbers $k$ and $k+1$ lie on the same upslope in $w$. Then $u$ is obtained from $w$ by moving $k+1$ to the left. 

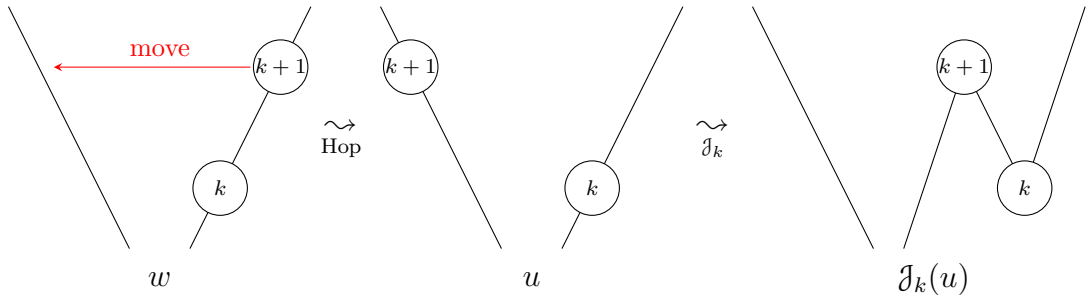
\begin{center}
\begin{tikzpicture}[scale=0.4]
\draw (-5,10)--(-1,2);
\draw (1,2)--(5,10);
\draw[fill=white] (2,4) circle (0.9cm);
\draw[fill=white] (4,8) circle (0.9cm);
\node at (2,4) {\tiny$k$};
\node at (4,8) {\tiny$k+1$};
\draw[-stealth, red] (3,8)--(-3.5,8);
\node[above] at (0,8) {\footnotesize\red{move}};
\node at (0,1) {$w$};

\node[below] at (6,6) {\tiny{Hop}};
\node at (6,6) {$\leadsto$};

\begin{scope}[xshift=350]
\draw (-5,10)--(-1,2);
\draw (1,2)--(5,10);
\draw[fill=white] (2,4) circle (0.9cm);
\node at (2,4) {\tiny$k$};

\draw[fill=white] (-4,8) circle (0.9cm);
\node at (-4,8) {\tiny$k+1$};
\node at (0,1) {$u$};

\node[below] at (6,6) {\tiny$\mathcal{J}_k$};
\node at (6,6) {$\leadsto$};
\end{scope}

\begin{scope}[xshift=700]
\draw (-5,10)--(-1,2);
\draw (0,2)--(2,8)--(4,4)--(6,10);
\draw[fill=white] (4,4) circle (0.9cm);
\node at (4,4) {\tiny$k$};

\draw[fill=white] (2,8) circle (0.9cm);
\node at (2,8) {\tiny$k+1$};
\node at (1,1) {$\mathcal{J}_k(u)$};
\end{scope}

\end{tikzpicture}

\end{center}

\noindent
\textbf{Case (b):} The number $k$ is a peak of $w$ and $k+1$ lies on an upslope. This implies that there are no mountains with the height greater than or equal to $k+2$ between the upslope containing $k+1$ and the mountain whose peak is $k$. Then $u$ is obtained from $w$ by moving $k+1$ to the left.

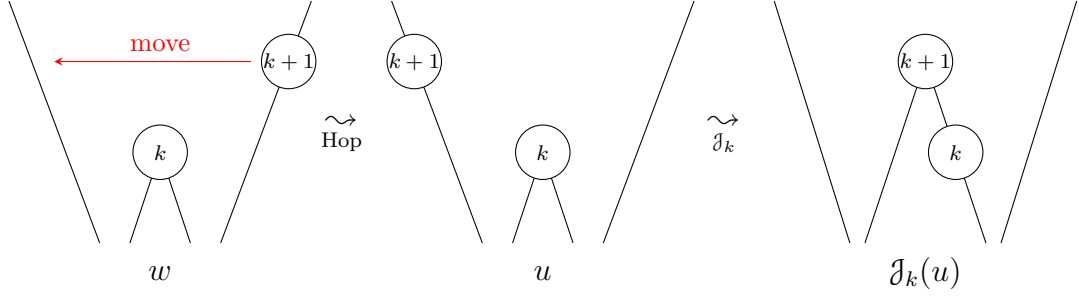
\begin{center}
\begin{tikzpicture}[scale=0.4]
\draw (-5,10)--(-2,2);
\draw (2,2)--(5,10);
\draw (-1, 2)--(0,5)--(1,2);
\draw[fill=white] (0,5) circle (0.9cm);
\draw[fill=white] (2+9/4,8) circle (0.9cm);
\node at (0,5) {\tiny$k$};
\node at (2+9/4,8) {\tiny$k+1$};
\draw[-stealth, red] (3,8)--(-3.5,8);
\node[above] at (0,8) {\footnotesize\red{move}};
\node at (0,1) {$w$};

\node[below] at (6,6) {\tiny{Hop}};
\node at (6,6) {$\leadsto$};

\begin{scope}[xshift=360]
\draw (-5,10)--(-2,2);
\draw (2,2)--(5,10);
\draw (-1, 2)--(0,5)--(1,2);
\draw[fill=white] (0,5) circle (0.9cm);
\draw[fill=white] (-2-9/4,8) circle (0.9cm);
\node at (0,5) {\tiny$k$};
\node at (-2-9/4,8) {\tiny$k+1$};

\node at (0,1) {$u$};

\node[below] at (6,6) {\tiny$\mathcal{J}_k$};
\node at (6,6) {$\leadsto$};
\end{scope}

\begin{scope}[xshift=720]
\draw (-5,10)--(-2.5,2);
\draw (2.5,2)--(5,10);
\draw (-2, 2)--(0,8)--(2,2);
\draw[fill=white] (1,5) circle (0.9cm);
\draw[fill=white] (0,8) circle (0.9cm);
\node at (1,5) {\tiny$k$};
\node at (0,8) {\tiny$k+1$};

\node at (0,1) {$\mathcal{J}_k(u)$};
\end{scope}
\end{tikzpicture}
\end{center}

\smallskip
In both cases (a) and (b), we have 
\begin{align*}
\des(\J_k(u))=\des(u).
\end{align*}
Indeed, the sum of the number of peaks of $u$ and the number of free letters on downslopes of $u$ equals the same sum with respect to $\J_k(u)$.
This proves the base case $\ell=1$. 

Now suppose that $\ell > 1$ and assume by induction that the claim is true for
arbitrary closed intervals $K'$ of successive integers with $|K'| \leq \ell-1$.
If 
\[
\J_{k+i}\J_{k+i-1} \cdots \J_k(u) = \J_{k+i-1}\J_{k+i-2} \cdots \J_k(u)
\] 
for some $0 \leq i \leq \ell-1$, then  we have
\begin{align} \label{eq:proof descent-preserving for closed interval 1}
\J_K(u) = \J_{K''}\J_{K'}(u). 
\end{align}
for $K'=[k, k+i-1]$ and $K''=[k+i+1,k+\ell-1]$. 
If subletters of $k+i+1,\ldots,k+\ell-1$ move in a construction of $u$ from $w$ (by crossing valleys), then we construct $w'$ from $\J_{K'}(u)$ by moving the subletters. 
One can see that $w' \in \widetilde{W^{K''}}$ with $\J_{K'}(u) \in \Hop(w')$. 
By the inductive hypothesis, we obtain 
\begin{align} \label{eq:proof descent-preserving for closed interval 2}
\des(\J_{K''}\J_{K'}(u)) = \des(\J_{K'}(u)). 
\end{align}
Noting that $w \in \widetilde{W^{K}} \subset \widetilde{W^{K'}}$, it also follows from the inductive assumption that 
\begin{align} \label{eq:proof descent-preserving for closed interval 3}
\des(\J_{K'}(u)) = \des(u). 
\end{align}
By \eqref{eq:proof descent-preserving for closed interval 1}, \eqref{eq:proof descent-preserving for closed interval 2} and \eqref{eq:proof descent-preserving for closed interval 3}, we have the desired equality $\des(\J_K(u))=\des(u)$.

Next, we assume that 
\begin{align*}
\J_{k+i}\J_{k+i-1} \cdots \J_k(u) \neq \J_{k+i-1}\J_{k+i-2} \cdots \J_k(u) \ \ \ {\rm for \ all} \ 0 \leq i \leq \ell-1,
\end{align*}
where we take by convention $\J_{k+i-1}\J_{k+i-2} \cdots \J_k=id$ whenever $i=0$.
This situation only happens in the following two cases.

\smallskip 
\noindent
\textbf{Case (i):} The numbers $k,k+1,\dots ,k+\ell$ lie on the same upslope of $w$, and $u$ is obtained from $w$ by moving $k+1,k+2,\dots ,k+\ell$ to the left.  
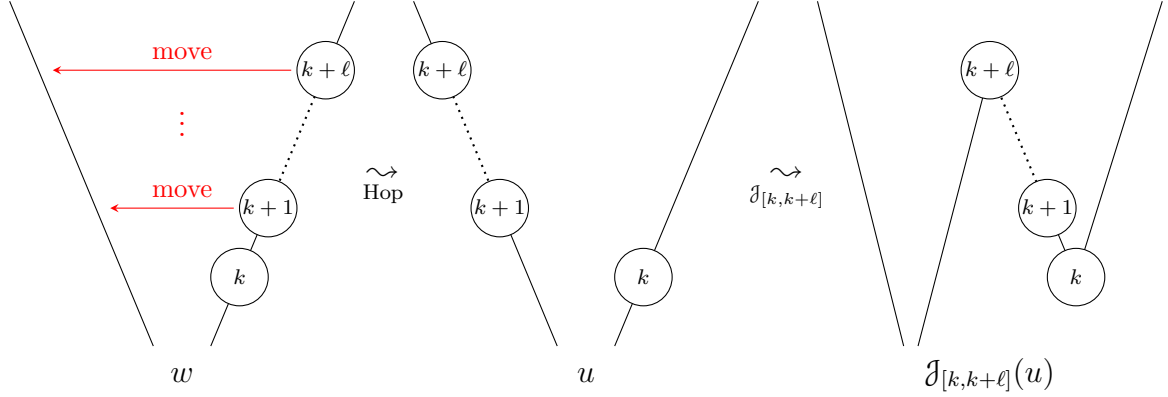
\begin{center}
\begin{tikzpicture}[scale=0.38]
\draw (-6,12)--(-1,0);
\draw (1,0)--(3,24/5);
\draw[dotted, thick] (3,24/5)--(5,48/5);
\draw (5,48/5)--(6,12);

\draw[fill=white] (2,12/5) circle (1cm);
\draw[fill=white] (3,24/5) circle (1cm);
\draw[fill=white] (5,48/5) circle (1cm);
\node at (2,12/5) {\tiny$k$};
\node at (3,24/5) {\tiny$k+ 1$};
\node at (5,48/5) {\tiny$k+ \ell$};

\draw[-stealth, red] (1.8,24/5)--(-2.5,24/5);
\node[above] at (0,24/5) {\footnotesize\red{move}};

\draw[-stealth, red] (3.8,48/5)--(-4.5,48/5);
\node[above] at (0,48/5) {\footnotesize\red{move}};

\node at (0,8) {\red{$\vdots$}};
\node at (0,-1) {$w$};

\node[below] at (7,6) {\tiny{Hop}};
\node at (7,6) {$\leadsto$};

\begin{scope}[xshift=400]
\draw (6,12)--(1,0);
\draw (-1,0)--(-3,24/5);
\draw[dotted, thick] (-3,24/5)--(-5,48/5);
\draw (-5,48/5)--(-6,12);

\draw[fill=white] (2,12/5) circle (1cm);
\draw[fill=white] (-3,24/5) circle (1cm);
\draw[fill=white] (-5,48/5) circle (1cm);
\node at (2,12/5) {\tiny$k$};
\node at (-3,24/5) {\tiny$k + 1$};
\node at (-5,48/5) {\tiny$k +\ell$};
\node at (0,-1) {$u$};

\node[below] at (7,6) {\tiny{$\mathcal{J}_{[k, k+\ell]}$}};
\node at (7,6) {$\leadsto$};
\end{scope}

\begin{scope}[xshift=800]
\draw (-6,12)--(-3,0);
\draw (-2.5,0)--(0,48/5);
\draw[thick, dotted] (0,48/5)--(2,24/5);
\draw (2,24/5)--(3,12/5)--(6,12);

\draw[fill=white] (3,12/5) circle (1cm);
\node at (3,12/5) {\tiny$k$};

\draw[fill=white] (2,24/5) circle (1cm);
\node at (2,24/5) {\tiny$k +1$};

\draw[fill=white] (-0,48/5) circle (1cm);
\node at (-0,48/5) {\tiny$k + \ell$};

\node at (0,-1) {$\mathcal{J}_{[k, k+\ell]}(u)$};
\end{scope}
\end{tikzpicture}
\end{center}

\noindent
\textbf{Case (ii):} The number $k$ is a peak of $w$ and $k+1,\dots, k+\ell$ lie on the same upslope with no mountains of height more than $k+\ell+1$ between the upslope containing $k+1, \dots, k+\ell$ and the mountain with peak $k$.  Moreover, $u$ is obtained from $w$ by moving $k+1,k+2,\dots, k+\ell$ to the left.

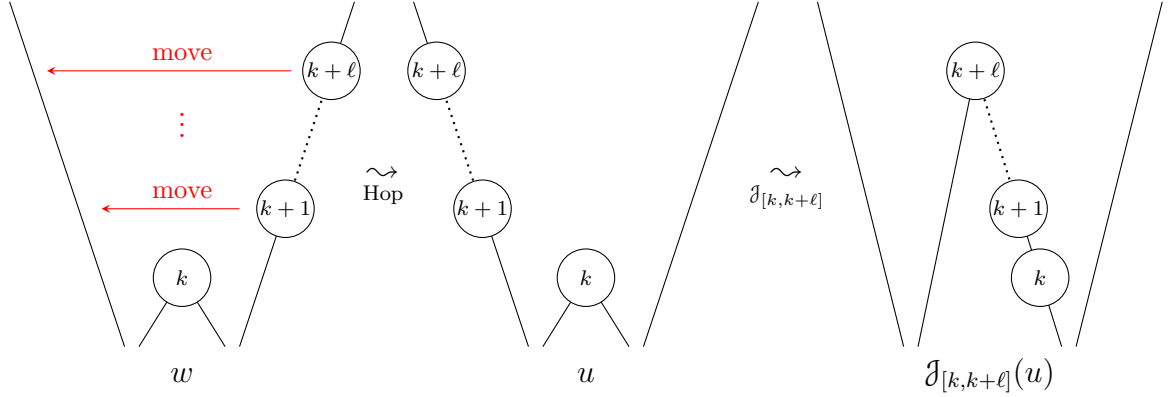
\begin{center}
\begin{tikzpicture}[scale=0.38]
\draw (-6,12)--(-2,0);
\draw (2,0)--(2+8/5,24/5); 
\draw[dotted, thick] (2+8/5,24/5)--(2+16/5,48/5);
\draw (2+16/5,48/5)--(6,12);
\draw (-1.5,0)--(0,12/5)--(1.5,0);

\draw[fill=white] (0,12/5) circle (1cm);
\node at (0,12/5) {\tiny$k$};

\draw[fill=white] (2+8/5,24/5) circle (1cm);
\node at (2+8/5,24/5) {\tiny$k+ 1$};

\draw[fill=white] (2+16/5,48/5) circle (1cm);
\node at (2+16/5,48/5) {\tiny$k+\ell$};

\node at (0,-1) {$w$};
\node at (7,6) {$\leadsto$};
\node[below] at (7,6) {\tiny{Hop}};

\draw[-stealth, red] (2,24/5)--(-2.8,24/5);
\node[above] at (0,24/5) {\footnotesize\red{move}};
\node at (0,8) {\red{$\vdots$}};
\draw[-stealth, red] (3.8,48/5)--(-4.7,48/5);
\node[above] at (0,48/5) {\footnotesize\red{move}};

\begin{scope}[xshift=400]
\draw (6,12)--(2,0);
\draw (-2,0)--(-2-8/5,24/5); 
\draw[dotted, thick] (-2-8/5,24/5)--(-2-16/5,48/5);
\draw (-2-16/5,48/5)--(-6,12);
\draw (-1.5,0)--(0,12/5)--(1.5,0);

\draw[fill=white] (0,12/5) circle (1cm);
\node at (0,12/5) {\tiny$k$};

\draw[fill=white] (-2-8/5,24/5) circle (1cm);
\node at (-2-8/5,24/5) {\tiny$k+1$};

\draw[fill=white] (-2-16/5,48/5) circle (1cm);
\node at (-2-16/5,48/5) {\tiny$k+ \ell$};

\node at (0,-1) {$u$};
\node at (7,6) {$\leadsto$};
\node[below] at (7,6) {\tiny{$\mathcal{J}_{[k, k+\ell]}$}};
\end{scope}

\begin{scope}[xshift=800]
\draw (-6,12)--(-3,0);
\draw (-2.5,0)--(-1/2,48/5);
\draw[thick, dotted] (-1/2,48/5)--(1,24/5);
\draw (1,24/5)--(7/4,12/5)--(2.5,0);
\draw (6,12)--(3,0);

\draw[fill=white] (7/4,12/5) circle (1cm);
\node at (7/4,12/5) {\tiny$k$};

\draw[fill=white] (1,24/5) circle (1cm);
\node at (1,24/5) {\tiny$k+ 1$};

\draw[fill=white] (-1/2,48/5) circle (1cm);
\node at (-1/2,48/5) {\tiny$k+ \ell$};

\node at (0,-1) {$\mathcal{J}_{[k, k+\ell]}(u)$};
\end{scope}
\end{tikzpicture}
\end{center}
In both cases (i) and (ii), one can see from \eqref{eq:descent} that 
\[
\des (\J_{k+i}\J_{k+i-1} \cdots \J_k(u))=\des(\J_{k+i-1}\J_{k+i-2} \cdots \J_k(u))
\] 
for all $0 \leq i \leq \ell-1$. This implies that $\des(\J_K(u))=\des(u)$, as desired. 
\end{proof}

In order to prove the claim $\des (\J_K(u))=\des (u)$ for arbitrary $K$, we prepare with the following definition.

\begin{definition}\label{eq_JK_operation}
Let $K$ be a subset of $[n-1]$.
For a permutation $u \in \mathfrak{S}_n$ with $\J_K(u) \neq u$,
one can write 
\begin{align} \label{eq:effective expression1}
\J_K(u) = \J_{i_p} \cdots \J_{i_2} \J_{i_1} (u)
\end{align}
for some $\{i_1, i_2, \dots, i_p \} \subseteq K$ with $i_1<i_2<\cdots<i_p$ such that 
\[
\J_{i_a} \big(\J_{i_{a-1}} \cdots \J_{i_1} (u)\big) \neq \J_{i_{a-1}} \cdots \J_{i_1} (u)
\]
for any $1 \leq a \leq p$. 
(Here, we take the convention $\J_{i_1} (u) \neq u$ whenever $a=1$.)
Then, we can rewrite \eqref{eq:effective expression1} as  
\begin{align} \label{eq:effective expression2}
\J_K(u)=\J_{K_m} \cdots \J_{K_1}(u)
\end{align}
for some closed intervals $K_j \coloneqq [k_j,k_j+\ell_j-1]$ in $K$ so that all of 
\begin{itemize} 
\item $k_i<k_i+\ell_i<k_{i+1}$ for $i \in [m-1]$ and $\ell_m>0$; 
\item  $\J_{k_j+i}\J_{k_j+i-1} \cdots \J_{k_j}(u_{j}) \neq \J_{k_j+i-1}\J_{k_j+i-2} \cdots \J_{k_j}(u_{j})$ for all $0 \leq i \leq \ell_j-1$ with $u_1\coloneqq u$ and $u_{j} \coloneqq \J_{K_{j-1}} \cdots \J_{K_1}(u)$ for $2 \leq j \leq m$.
\end{itemize} 
We call the expression in \eqref{eq:effective expression2} the \emph{effective expression for $\J_K(u)$}. 
Note that closed intervals $K_j=[k_j,k_j+\ell_j-1]$ in $K$ depend on $u$ and $K_j$ is not necessary a maximal consecutive string of $K$.
\end{definition} 

\begin{example}
Let $n=9$ and $K=\{2,3,4,6,7,8 \}$. For $w=672813459 \in \twk$ and $u=H_9(w)= 967284135 \in \Hop(w)$ (as shown in Figure~\ref{fig_valley_hopping}), one can see that 
$$
\J_K(u) =672981435 =  \J_8 \J_3(u)
$$
where $\J_3(u) \neq u$ and $\J_8 \J_3(u) \neq \J_3(u)$.
Also, it is the effective expression for $\J_K(u)$.
\end{example}

\begin{proposition} \label{proposition:descent-preserving}
Let $K$ be an arbitrary subset of $[n-1]$.
Let $w \in \twk$ and $u \in \Hop(w)$. 
Then we have
\begin{align*}
\des(\J_K(u))=\des(u).
\end{align*}
\end{proposition}

\begin{proof}
We may assume that $\J_K(u) \neq u$.
Consider the effective expression for $\J_K(u)$ in \eqref{eq:effective expression2}.
By a similar argument to that in the proof of Lemma~\ref{lemma:descent-preserving for closed interval}, one can see that either Case (i) or Case (ii) below happens in the mountain range description of $u \in \Hop(w)$. 
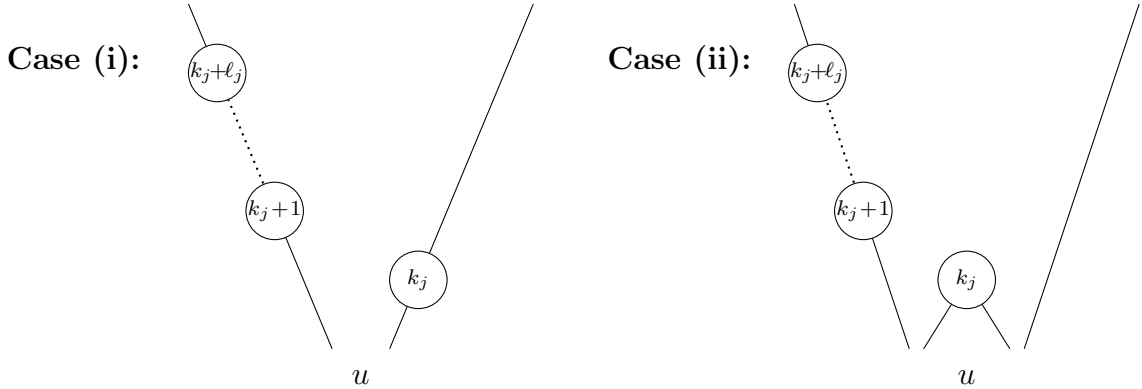
\[
\begin{tikzpicture}[scale=0.38, baseline=(current bounding box.north)]
\node at (-10, 10) {\textbf{Case (i): }};
\draw (6,12)--(1,0);
\draw (-1,0)--(-3,24/5);
\draw[dotted, thick] (-3,24/5)--(-5,48/5);
\draw (-5,48/5)--(-6,12);

\draw[fill=white] (2,12/5) circle (1cm);
\draw[fill=white] (-3,24/5) circle (1cm);
\draw[fill=white] (-5,48/5) circle (1cm);
\node at (2,12/5) {\tiny$k_j$};
\node at (-3,24/5) {\tiny$k_j \!  + \! 1$};
\node at (-5,48/5) {\tiny$k_j \! \! + \! \! \ell_j$};
\node at (0,-1) {$u$};
\begin{scope}[xshift=600]
\node at (-10, 10) {\textbf{Case (ii): }};
\draw (6,12)--(2,0);
\draw (-2,0)--(-2-8/5,24/5); 
\draw[dotted, thick] (-2-8/5,24/5)--(-2-16/5,48/5);
\draw (-2-16/5,48/5)--(-6,12);
\draw (-1.5,0)--(0,12/5)--(1.5,0);

\draw[fill=white] (0,12/5) circle (1cm);
\node at (0,12/5) {\tiny$k_j$};

\draw[fill=white] (-2-8/5,24/5) circle (1cm);
\node at (-2-8/5,24/5) {\tiny$k_j \! +  \! 1$};

\draw[fill=white] (-2-16/5,48/5) circle (1cm);
\node at (-2-16/5,48/5) {\tiny$k_j \! \! + \! \! \ell_j$};

\node at (0,-1) {$u$};
\end{scope}
\end{tikzpicture}
\]

%
%
%
%
%
%
%
%
%
%
However, by the definition of $u_j$, the positions of the letters $k_j,k_j+1,\ldots,k_j+\ell_j$ in the mountain range description of $u_j$ do not change from their positions in $u$ for each $1 \leq j \leq m$. 
Note that the letters $k_j+1,\ldots,k_j+\ell_j$ in the mountain range description of $u_j$ are free letters on the same downslope.
Thus, we can find the element $w_j \in \widetilde{W^{K_j}}$ such that $u_j \in \Hop(w_j)$ which is obtained from $u_j$ by moving all free letters which lie on a downslope (by crossing valleys). (See Figures~\ref{picture:wj Case(i)} and \ref{picture:wj Case(ii)}.)
%
%
%
%
%
%
%
%
%

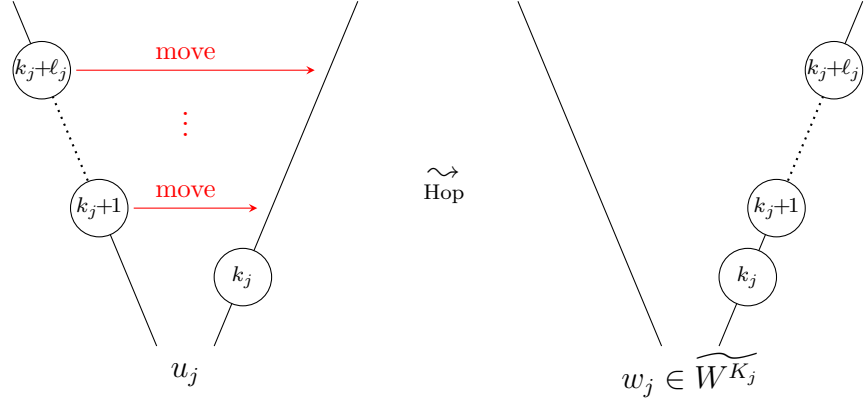
\begin{figure}
\begin{tikzpicture}[scale=0.38]
\draw (6,12)--(1,0);
\draw (-1,0)--(-3,24/5);
\draw[dotted, thick] (-3,24/5)--(-5,48/5);
\draw (-5,48/5)--(-6,12);

\draw[fill=white] (2,12/5) circle (1cm);
\draw[fill=white] (-3,24/5) circle (1cm);
\draw[fill=white] (-5,48/5) circle (1cm);
\node at (2,12/5) {\tiny$k_j$};
\node at (-3,24/5) {\tiny$k_j \! \! + \! \! 1$};
\node at (-5,48/5) {\tiny$k_j \! \! + \! \! \ell_j$};

\draw[-stealth, red](-1.8,24/5)--(2.5,24/5);
\node[above] at (0,24/5) {\footnotesize\red{move}};

\draw[-stealth, red] (-3.8,48/5)--(4.5,48/5);
\node[above] at (0,48/5) {\footnotesize\red{move}};

\node at (0,8) {\red{$\vdots$}};
\node at (0,-1) {$u_j$};

\node[below] at (9,6) {\tiny{Hop}};
\node at (9,6) {$\leadsto$};

\begin{scope}[xshift=500]
\draw (-6,12)--(-1,0);
\draw (1,0)--(3,24/5);
\draw[dotted, thick] (3,24/5)--(5,48/5);
\draw (5,48/5)--(6,12);

\draw[fill=white] (2,12/5) circle (1cm);
\draw[fill=white] (3,24/5) circle (1cm);
\draw[fill=white] (5,48/5) circle (1cm);
\node at (2,12/5) {\tiny$k_j$};
\node at (3,24/5) {\tiny$k_j \! \! + \! \! 1$};
\node at (5,48/5) {\tiny$k_j \! \! + \! \! \ell_j$};

\node at (0,-1) {$w_j\in \widetilde{W^{K_j}}$};
\end{scope}
\end{tikzpicture}
\caption{$w_j \in \widetilde{W^{K_j}}$ in Case (i).}
\label{picture:wj Case(i)}
\end{figure}

\begin{figure}[h]
\begin{tikzpicture}[scale=0.38]
\draw (6,12)--(2,0);
\draw (-2,0)--(-2-8/5,24/5); 
\draw[dotted, thick] (-2-8/5,24/5)--(-2-16/5,48/5);
\draw (-2-16/5,48/5)--(-6,12);
\draw (-1.5,0)--(0,12/5)--(1.5,0);

\draw[fill=white] (0,12/5) circle (1cm);
\node at (0,12/5) {\tiny$k_j$};

\draw[fill=white] (-2-8/5,24/5) circle (1cm);
\node at (-2-8/5,24/5) {\tiny$k_j \! \! + \! \! 1$};

\draw[fill=white] (-2-16/5,48/5) circle (1cm);
\node at (-2-16/5,48/5) {\tiny$k_j \! \! + \! \! \ell_j$};

\draw[-stealth, red] (-2,24/5)--(2.8,24/5);
\node[above] at (0,24/5) {\footnotesize\red{move}};
\node at (0,8) {\red{$\vdots$}};
\draw[-stealth, red] (-3.8,48/5)--(4.7,48/5);
\node[above] at (0,48/5) {\footnotesize\red{move}};

\node at (0,-1) {$u_j$};

\node[below] at (9,6) {\tiny{Hop}};
\node at (9,6) {$\leadsto$};

\begin{scope}[xshift=500]
\draw (-6,12)--(-2,0);
\draw (2,0)--(2+8/5,24/5); 
\draw[dotted, thick] (2+8/5,24/5)--(2+16/5,48/5);
\draw (2+16/5,48/5)--(6,12);
\draw (-1.5,0)--(0,12/5)--(1.5,0);

\draw[fill=white] (0,12/5) circle (1cm);
\node at (0,12/5) {\tiny$k_j$};

\draw[fill=white] (2+8/5,24/5) circle (1cm);
\node at (2+8/5,24/5) {\tiny$k_j \! \! + \! \! 1$};

\draw[fill=white] (2+16/5,48/5) circle (1cm);
\node at (2+16/5,48/5) {\tiny$k_j \! \! + \! \! \ell_j$};

\node at (0,-1) {$w_j\in \widetilde{W^{K_j}}$};
\end{scope}

%
%
%
%
%
%
%
%
%
%
\end{tikzpicture}
\caption{$w_j \in \widetilde{W^{K_j}}$ in Case (ii).}
\label{picture:wj Case(ii)}
\end{figure}
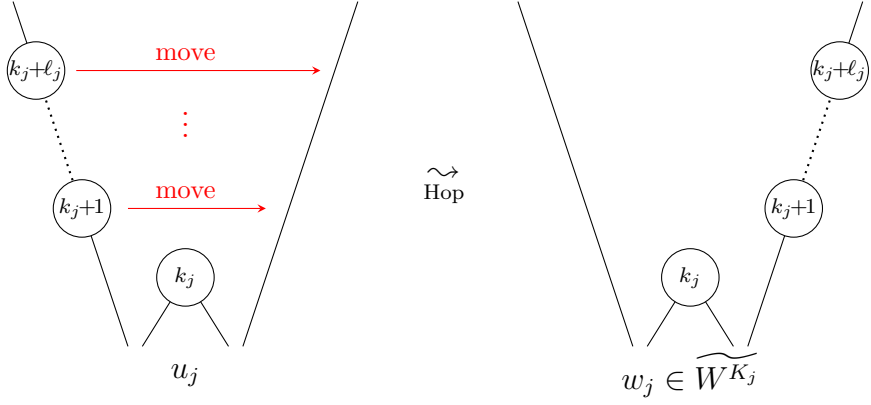

Now, we apply Lemma \ref{lemma:descent-preserving for closed interval} repeatedly to 
conclude  
\[
\des(\J_K(u)) =\des(\J_{K_{m}}(u_{m}))=\des(u_{m})=\des(\J_{K_{m-1}}(u_{m-1}))= \cdots = \des(u_1)=\des(u).
\]
\end{proof}

\section{Inverse operator}\label{sec_Inv_oper}
In this section we prove the bijectivity of the function $\Theta$ defined in Proposition \ref{proposition:key} by constructing its inverse map. Let $K$ be a subset of $[n-1]$. Take permutations $w\in \twk$ and $u\in \Hop(w)$ such that $\J_K(u)\neq u$. Consider the effective expression for $\J_K(u)$ defined in \eqref{eq:effective expression2}.
Then, for each $j$, we have the following two cases for the positions of the letters $k_j,k_j+1,\ldots,k_j+\ell_j$ in the mountain range descriptions for $w$, $u$ and $\J_K(u)$. 

\bigskip 
\noindent
\textbf{Case (i):} 

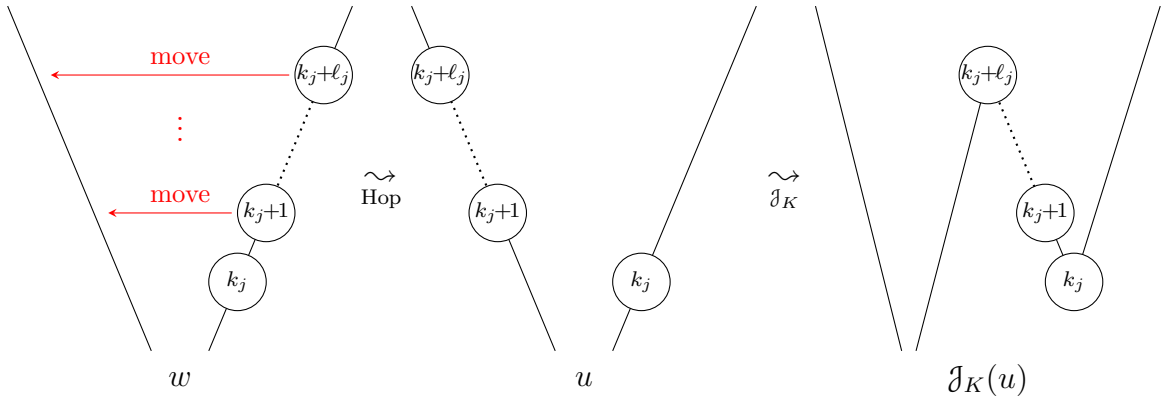
\begin{center}
\begin{tikzpicture}[scale=0.38]
\draw (-6,12)--(-1,0);
\draw (1,0)--(3,24/5);
\draw[dotted, thick] (3,24/5)--(5,48/5);
\draw (5,48/5)--(6,12);

\draw[fill=white] (2,12/5) circle (1cm);
\draw[fill=white] (3,24/5) circle (1cm);
\draw[fill=white] (5,48/5) circle (1cm);
\node at (2,12/5) {\tiny$k_j$};
\node at (3,24/5) {\tiny$k_j \! \! + \! \! 1$};
\node at (5,48/5) {\tiny$k_j \! \! + \! \! \ell_j$};

\draw[-stealth, red] (1.8,24/5)--(-2.5,24/5);
\node[above] at (0,24/5) {\footnotesize\red{move}};

\draw[-stealth, red] (3.8,48/5)--(-4.5,48/5);
\node[above] at (0,48/5) {\footnotesize\red{move}};

\node at (0,8) {\red{$\vdots$}};
\node at (0,-1) {$w$};

\node[below] at (7,6) {\tiny{Hop}};
\node at (7,6) {$\leadsto$};

\begin{scope}[xshift=400]
\draw (6,12)--(1,0);
\draw (-1,0)--(-3,24/5);
\draw[dotted, thick] (-3,24/5)--(-5,48/5);
\draw (-5,48/5)--(-6,12);

\draw[fill=white] (2,12/5) circle (1cm);
\draw[fill=white] (-3,24/5) circle (1cm);
\draw[fill=white] (-5,48/5) circle (1cm);
\node at (2,12/5) {\tiny$k_j$};
\node at (-3,24/5) {\tiny$k_j \! \! + \! \! 1$};
\node at (-5,48/5) {\tiny$k_j \! \! + \! \! \ell_j$};
\node at (0,-1) {$u$};

\node[below] at (7,6) {\tiny{$\mathcal{J}_K$}};
\node at (7,6) {$\leadsto$};
\end{scope}

\begin{scope}[xshift=800]
\draw (-6,12)--(-3,0);
\draw (-2.5,0)--(0,48/5);
\draw[thick, dotted] (0,48/5)--(2,24/5);
\draw (2,24/5)--(3,12/5)--(6,12);

\draw[fill=white] (3,12/5) circle (1cm);
\node at (3,12/5) {\tiny$k_j$};

\draw[fill=white] (2,24/5) circle (1cm);
\node at (2,24/5) {\tiny$k_j \! \! + \! \! 1$};

\draw[fill=white] (-0,48/5) circle (1cm);
\node at (-0,48/5) {\tiny$k_j \! \! + \! \! \ell_j$};

\node at (0,-1) {$\mathcal{J}_K(u)$};
\end{scope}
\end{tikzpicture}
\end{center}

%
%
%
%
%
%
%
%
%
%
%
%
%
%
%
%


\noindent
\textbf{Case (ii):} 
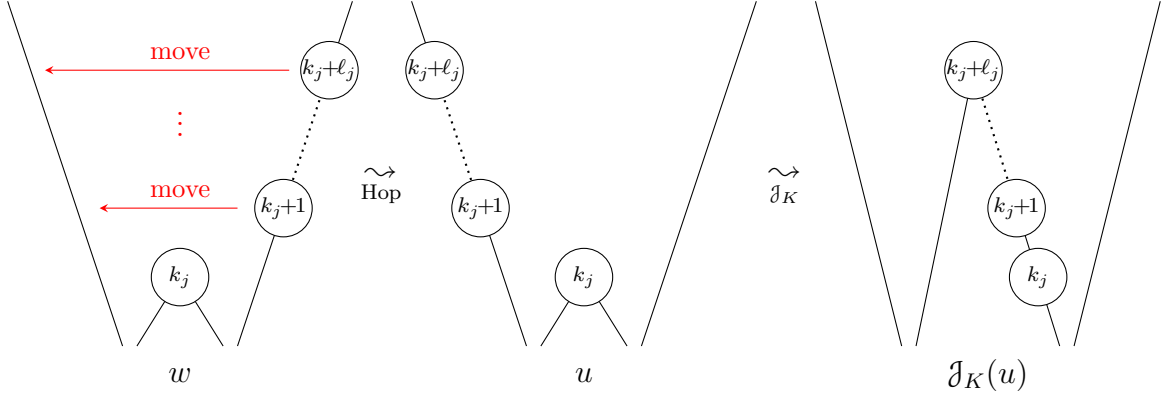
\begin{center}
\begin{tikzpicture}[scale=0.38]
\draw (-6,12)--(-2,0);
\draw (2,0)--(2+8/5,24/5); 
\draw[dotted, thick] (2+8/5,24/5)--(2+16/5,48/5);
\draw (2+16/5,48/5)--(6,12);
\draw (-1.5,0)--(0,12/5)--(1.5,0);

\draw[fill=white] (0,12/5) circle (1cm);
\node at (0,12/5) {\tiny$k_j$};

\draw[fill=white] (2+8/5,24/5) circle (1cm);
\node at (2+8/5,24/5) {\tiny$k_j \! \! + \! \! 1$};

\draw[fill=white] (2+16/5,48/5) circle (1cm);
\node at (2+16/5,48/5) {\tiny$k_j \! \! + \! \! \ell_j$};

\node at (0,-1) {$w$};
\node at (7,6) {$\leadsto$};
\node[below] at (7,6) {\tiny{Hop}};

\draw[-stealth, red] (2,24/5)--(-2.8,24/5);
\node[above] at (0,24/5) {\footnotesize\red{move}};
\node at (0,8) {\red{$\vdots$}};
\draw[-stealth, red] (3.8,48/5)--(-4.7,48/5);
\node[above] at (0,48/5) {\footnotesize\red{move}};

\begin{scope}[xshift=400]
\draw (6,12)--(2,0);
\draw (-2,0)--(-2-8/5,24/5); 
\draw[dotted, thick] (-2-8/5,24/5)--(-2-16/5,48/5);
\draw (-2-16/5,48/5)--(-6,12);
\draw (-1.5,0)--(0,12/5)--(1.5,0);

\draw[fill=white] (0,12/5) circle (1cm);
\node at (0,12/5) {\tiny$k_j$};

\draw[fill=white] (-2-8/5,24/5) circle (1cm);
\node at (-2-8/5,24/5) {\tiny$k_j \! \! + \! \! 1$};

\draw[fill=white] (-2-16/5,48/5) circle (1cm);
\node at (-2-16/5,48/5) {\tiny$k_j \! \! + \! \! \ell_j$};

\node at (0,-1) {$u$};
\node at (7,6) {$\leadsto$};
\node[below] at (7,6) {\tiny{$\mathcal{J}_K$}};
\end{scope}

\begin{scope}[xshift=800]
\draw (-6,12)--(-3,0);
\draw (-2.5,0)--(-1/2,48/5);
\draw[thick, dotted] (-1/2,48/5)--(1,24/5);
\draw (1,24/5)--(7/4,12/5)--(2.5,0);
\draw (6,12)--(3,0);

\draw[fill=white] (7/4,12/5) circle (1cm);
\node at (7/4,12/5) {\tiny$k_j$};

\draw[fill=white] (1,24/5) circle (1cm);
\node at (1,24/5) {\tiny$k_j \! \! + \! \! 1$};

\draw[fill=white] (-1/2,48/5) circle (1cm);
\node at (-1/2,48/5) {\tiny$k_j \! \! + \! \! \ell_j$};

\node at (0,-1) {$\mathcal{J}_K(u)$};
\end{scope}
\end{tikzpicture}
\end{center}

Roughly speaking, the operator $\J_K$ means either ``building a new mountain'' (Case~(i)), or  ``growing a height of a mountain'' (Case~(ii)).
We first have the following lemma. 

\begin{lemma} \label{lemma:kj-1}
Let $K, w$ and $u$ be as above. If  $(\J_K(u))^{-1}(k_j-1)=(\J_K(u))^{-1}(k_j)+1$ (which happens only in Case (ii) in the first paragraph of this section), then $k_j-1 \not\in K$.
\end{lemma}

\begin{proof}
If $k_j-1 \in K$, then we have $u^{-1}(k_j-1) = u^{-1}(k_j)+1$.
Indeed, if we write $\J_K(u)=\J_{K_m}\cdots \J_{K_1}(u)$ as the effective expression discussed in \eqref{eq:effective expression2}, then the assumption $(\J_K(u))^{-1}(k_j-1)=(\J_K(u))^{-1}(k_j)+1$ implies that $k_j-1$ does not belong to any closed intervals $K_i's$. Since $w \in \Hop(u)$, one has $w^{-1}(k_j) < w^{-1}(k_j-1)$. 
However, this is a contradiction since $k_j-1, k_j$ have to appear in the one-line notation of $w \in \twk$ as a subsequence in increasing order by \eqref{eq:KW}.
Hence, we have $k_j-1 \not\in K$.
\end{proof}

We observe that 
\begin{enumerate}
\item[$\bullet$] $k_j+\ell_j$ is a peak of $\J_K(u)$ in both Cases (i) and (ii).
\item[$\bullet$] A slope immediately after $k_j$ is an upslope in Case (i), while a slope immediately after $k_j$ is a downslope in Case (ii).
\end{enumerate}

Based on these two facts, we now construct an operator $\mathcal{L}_K$ on $W(K)$ with $\mathcal{L}_K(v)\in \mathfrak{S}_n$ for $v\in W(K)$. For $K\subset [n-1]$ and $v\in W(K)$, we denote by 
\begin{equation}\label{eq_string}
(k+\ell; k+\ell-1, \dots, k)
\end{equation}
a consecutive string in $K$ consisting of a peak $k+\ell$ of $v$ and letters $k+\ell-1, \dots, k$ immediately followed by $k+\ell$ such that $\{k+\ell-1, \dots, k\}\subset K$ and $k-1\notin K$ whenever $v^{-1}(k-1)-v^{-1}(k)=1$.
\begin{definition} \label{definition:key inverse operator}
For $v\in W(K)$, consider all strings $(k_j+\ell_j; k_j+\ell_j-1, \dots, k_j)$ in $v$ of the form \eqref{eq_string}. If $v$ has no such a string, then we define $\L_K(v)=v$. Otherwise, assuming  $1\leq j \leq m$ and $k_1 < k_1+\ell_1 < k_2 < k_2+\ell_2 < \cdots < k_m < k_m+\ell_m$ without loss of generality, we define $\mathcal{L}_K(v)$ inductively as follows. 
\begin{enumerate}
\item $v_{m+1} \coloneqq v$.
\item Given $v_{j+1}$ with $j\leq m$, let $c_j \coloneqq v_{j+1}^{-1}(k_j+\ell_j)$ and $a_j$ the largest integer such that $a_j<c_j$ and  $v_{j+1}(a_j) >  k_j+\ell_j$. We define $v_{j}$ by the following one-line notation:
\[
\hspace{-25pt}
v_{j+1}(1) \cdots v_{j+1}(a_j) \, k_j+\ell_j \ \ k_j+\ell_j-1 \ \cdots \ k_j+1 \, v_{j+1}(a_j+1) \cdots v_{j+1}(c_j-1) \ \ k_j \ \ v_{j+1}(c_j+\ell_j+1) \cdots v_{j+1}(n).
\]
\item We define $\L_K(v) \coloneqq v_1$.
\end{enumerate}
\end{definition}

\begin{figure}
\begin{tikzpicture}[scale=0.38]
\draw (-6,12)--(-3,0);
\draw (-2.5,0)--(-1/2,48/5);
\draw[thick, dotted] (-1/2,48/5)--(1,24/5);
\draw (1,24/5)--(7/4,12/5)--(2.5,0);
\draw (6,12)--(3,0);

\draw[fill=white] (7/4,12/5) circle (1cm);
\node at (7/4,12/5) {\tiny$k_j$};

\draw[fill=white] (1,24/5) circle (1cm);
\node at (1,24/5) {\tiny$k_j \! \! + \! \! 1$};

\draw[fill=white] (-1/2,48/5) circle (1cm);
\node at (-1/2,48/5) {\tiny$k_j \! \! + \! \! \ell_j$};

\node at (0,-1) {$v$};

\node at (8,6) {$\leadsto$};
\node[below] at (8,6) {\tiny{$\mathcal{L}_K$}};

\begin{scope}[xshift=500]
\draw (6,12)--(2,0);
\draw (-2,0)--(-2-8/5,24/5); 
\draw[dotted, thick] (-2-8/5,24/5)--(-2-16/5,48/5);
\draw (-2-16/5,48/5)--(-6,12);
\draw (-1.5,0)--(0,12/5)--(1.5,0);

\draw[fill=white] (0,12/5) circle (1cm);
\node at (0,12/5) {\tiny$k_j$};

\draw[fill=white] (-2-8/5,24/5) circle (1cm);
\node at (-2-8/5,24/5) {\tiny$k_j \! \! + \! \! 1$};

\draw[fill=white] (-2-16/5,48/5) circle (1cm);
\node at (-2-16/5,48/5) {\tiny$k_j \! \! + \! \! \ell_j$};

\node at (0,-1) {$\mathcal{L}_K(u)$};
\end{scope}

\end{tikzpicture}
%
%
%
%
%
%
%
%
\caption{The $k_j$ is a free letter, namely the slope immediately after $k_j$ is a downslope.}
\label{picture:definition LKv-2}
\end{figure}
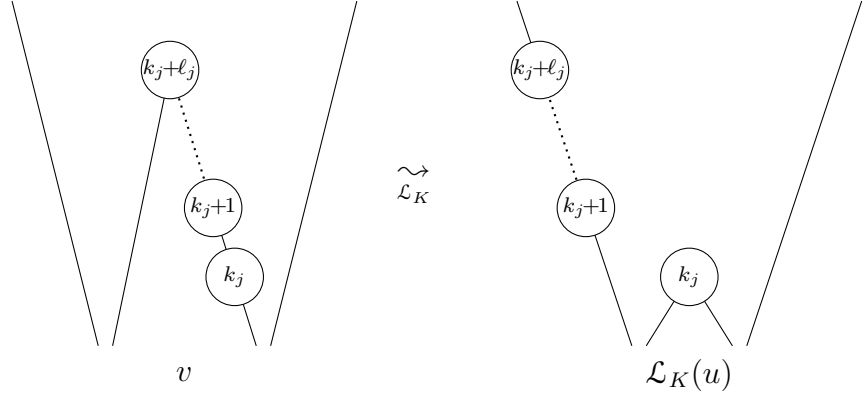

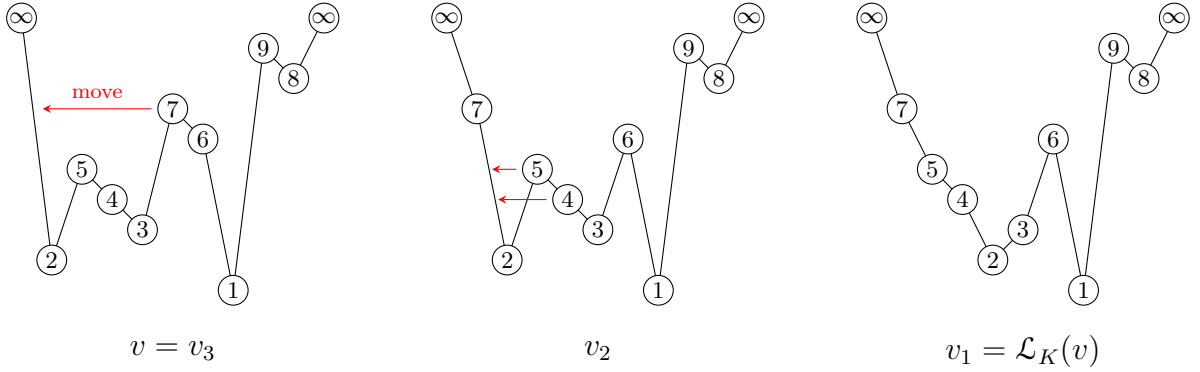
\begin{figure}
\begin{tikzpicture}[scale=0.4]
\draw (0,10)--(1,2)--(2,5)--(3,4)--(4,3)--(5,7)--(6,6)--(7,1)--(8,9)--(9,8)--(10,10);
\draw[fill=white] (0,10) circle (14pt); 
\draw[fill=white] (1,2) circle (14pt); 
\draw[fill=white] (2,5) circle (14pt); 
\draw[fill=white] (3,4) circle (14pt); 
\draw[fill=white] (4,3) circle (14pt); 
\draw[fill=white] (5,7) circle (14pt); 
\draw[fill=white] (6,6) circle (14pt); 
\draw[fill=white] (7,1) circle (14pt); 
\draw[fill=white] (8,9) circle (14pt); 
\draw[fill=white] (9,8) circle (14pt); 
\draw[fill=white] (10,10) circle (14pt); 

\node at (0,10) {\scriptsize$\infty$};
\node at (1,2) {\scriptsize$2$};
\node at (2,5) {\scriptsize$5$};
\node at (3,4) {\scriptsize$4$};
\node at (4,3) {\scriptsize$3$};
\node at (5,7) {\scriptsize$7$};
\node at (6,6) {\scriptsize$6$};
\node at (7,1) {\scriptsize$1$};
\node at (8,9) {\scriptsize$9$};
\node at (9,8) {\scriptsize$8$};
\node at (10,10) {\scriptsize$\infty$};

\draw[red, -stealth] (4.3,7)--(0.7,7);
\node[red, above] at (2.5,7) {\tiny{move}};

\node at (5,-1) {$v=v_3$};
\begin{scope}[xshift=400]
\draw (0,10)--(1,7)--(2,2)--(3,5)--(4,4)--(5,3)--(6,6)--(7,1)--(8,9)--(9,8)--(10,10);
\draw[fill=white] (0,10) circle (14pt); 
\draw[fill=white] (1,7) circle (14pt); 
\draw[fill=white] (2,2) circle (14pt); 
\draw[fill=white] (3,5) circle (14pt); 
\draw[fill=white] (4,4) circle (14pt); 
\draw[fill=white] (5,3) circle (14pt); 
\draw[fill=white] (6,6) circle (14pt); 
\draw[fill=white] (7,1) circle (14pt); 
\draw[fill=white] (8,9) circle (14pt); 
\draw[fill=white] (9,8) circle (14pt); 
\draw[fill=white] (10,10) circle (14pt); 

\node at (0,10) {\scriptsize$\infty$};
\node at (1,7) {\scriptsize$7$};
\node at (2,2) {\scriptsize$2$};
\node at (3,5) {\scriptsize$5$};
\node at (4,4) {\scriptsize$4$};
\node at (5,3) {\scriptsize$3$};
\node at (6,6) {\scriptsize$6$};
\node at (7,1) {\scriptsize$1$};
\node at (8,9) {\scriptsize$9$};
\node at (9,8) {\scriptsize$8$};
\node at (10,10) {\scriptsize$\infty$};

\draw[red, -stealth] (2.3,5)--(1.5,5);
\draw[red, -stealth] (3.3,4)--(1.7,4);

\node at (5,-1) {$v_2$};

\end{scope}

\begin{scope}[xshift=800]
\draw (0,10)--(1,7)--(2,5)--(3,4)--(4,2)--(5,3)--(6,6)--(7,1)--(8,9)--(9,8)--(10,10);
\draw[fill=white] (0,10) circle (14pt); 
\draw[fill=white] (1,7) circle (14pt); 
\draw[fill=white] (2,5) circle (14pt); 
\draw[fill=white] (3,4) circle (14pt); 
\draw[fill=white] (4,2) circle (14pt); 
\draw[fill=white] (5,3) circle (14pt); 
\draw[fill=white] (6,6) circle (14pt); 
\draw[fill=white] (7,1) circle (14pt); 
\draw[fill=white] (8,9) circle (14pt); 
\draw[fill=white] (9,8) circle (14pt); 
\draw[fill=white] (10,10) circle (14pt); 

\node at (0,10) {\scriptsize$\infty$};
\node at (1,7) {\scriptsize$7$};
\node at (2,5) {\scriptsize$5$};
\node at (3,4) {\scriptsize$4$};
\node at (4,2) {\scriptsize$2$};
\node at (5,3) {\scriptsize$3$};
\node at (6,6) {\scriptsize$6$};
\node at (7,1) {\scriptsize$1$};
\node at (8,9) {\scriptsize$9$};
\node at (9,8) {\scriptsize$8$};
\node at (10,10) {\scriptsize$\infty$};
\node at (5,-1) {$v_1=\L_K(v)$};
\end{scope}

\end{tikzpicture}
\caption{$\L_{\{3,4,6 \}}(254376198)$.}
\label{fig_LK254376198}
\end{figure}

We describe $\L_K(v)$ pictorially in Figures~\ref{picture:definition LKv-1} and \ref{picture:definition LKv-2}. 
Roughly speaking, the operator $\L_K$ means either ``leveling a hill'' (in Figure~\ref{picture:definition LKv-1}), or  ``reducing a height of a mountain'' (in Figure~\ref{picture:definition LKv-2}). For example,  if $K=\{3,4\} \sqcup \{6\} \subset [8]$ and $v=254376198\in W(K)\subset \mathfrak{S}_9$,  then we have $(5;4,3)$ and $(7;6)$ for strings of \eqref{eq_string}. This defines 
\[
v=v_3=254376198, ~ v_2=725436198, ~L_K(v)=v_1=754236198.
\] 
See Figure \ref{fig_LK254376198} for a pictorial description of this operation. On the other hand, applying the operation discussed in Definition \ref{eq_JK_operation}, we obtain $\J_{\{3,4\}}(v_1)=v_2$ and $\J_{\{6\}}(v_2)=v_3$, which means that $\J_K(\L_K(v))=v$. Indeed, we have the following proposition.

\begin{figure}
\begin{tikzpicture}[scale=0.38]
\draw (-6,12)--(-3,0);
\draw (-2.5,0)--(0,48/5);
\draw[thick, dotted] (0,48/5)--(2,24/5);
\draw (2,24/5)--(3,12/5)--(6,12);

\draw[fill=white] (3,12/5) circle (1cm);
\node at (3,12/5) {\tiny$k_j$};

\draw[fill=white] (2,24/5) circle (1cm);
\node at (2,24/5) {\tiny$k_j \! \! + \! \! 1$};

\draw[fill=white] (-0,48/5) circle (1cm);
\node at (-0,48/5) {\tiny$k_j \! \! + \! \! \ell_j$};

\node at (0,-1) {$v$};

\node[below] at (8,6) {\tiny{$\mathcal{L}_K$}};
\node at (8,6) {$\leadsto$};

\begin{scope}[xshift=500]
\draw (6,12)--(1,0);
\draw (-1,0)--(-3,24/5);
\draw[dotted, thick] (-3,24/5)--(-5,48/5);
\draw (-5,48/5)--(-6,12);

\draw[fill=white] (2,12/5) circle (1cm);
\draw[fill=white] (-3,24/5) circle (1cm);
\draw[fill=white] (-5,48/5) circle (1cm);
\node at (2,12/5) {\tiny$k_j$};
\node at (-3,24/5) {\tiny$k_j \! \! + \! \! 1$};
\node at (-5,48/5) {\tiny$k_j \! \! + \! \! \ell_j$};
\node at (0,-1) {$\mathcal{L}_K(v)$};
\end{scope}
\end{tikzpicture}
\caption{The letter $k_j$ is a valley, namely the slope immediately after $k_j$ is an upslope.}
\label{picture:definition LKv-1}
\end{figure}
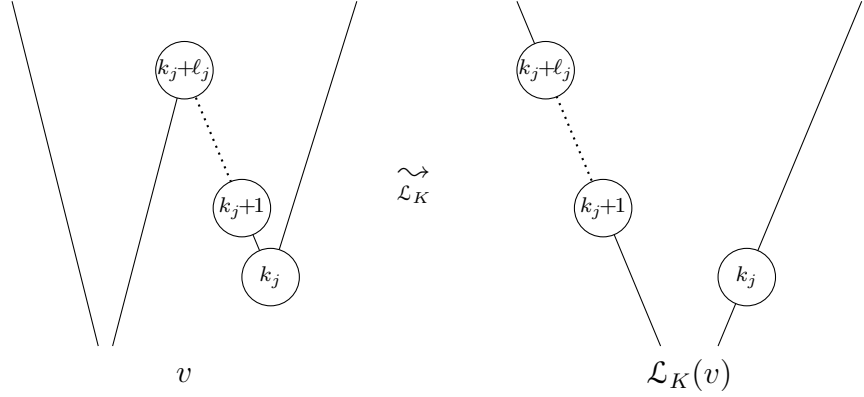

%
%
%
%
%
%
%
%

\begin{proposition} \label{proposition:bijection}
Let $K$ be a subset of $[n-1]$. 
\begin{enumerate}
\item \label{(1)} For $v \in W(K)$, the permutation $\L_K(v)\in\mathfrak{S}_n$ determines the unique permutation $w \in \twk$ such that $\L_K(v) \in \Hop(w)$.
\item The map
\begin{align*} 
W(K) \to \{(w,u) \in \twk \times \mathfrak{S}_n \mid u \in \Hop(w) \}
\end{align*}
defined by $v \mapsto (w,\L_K(v))$, where $w$ is the element determined by \eqref{(1)}, is the inverse of $\Theta$. 
\end{enumerate}
\end{proposition}

\begin{proof}
(1) Let $w$ be the canonical representative obtained from $\L_K(v)$ by moving free letters on a downslope to on an upslope by valley hopping. 
In other words, $w$ is the unique element of $\ts$ (see \eqref{eq:hatSn}) such that $\L_K(v) \in \Hop(w)$. 

It remains to prove that $w$ satisfies 
\eqref{eq:KW}.
Let $k$ be an arbitrary element of $K$.
It suffices to show that the numbers $k, k+1$ appear in the one-line notation of $w$ as a subsequence in increasing order. 
Since $v$ is an element of $W(K)$, either $k$ is the southwest of $k+1$ in the mountain range description  of $v$ or $k$ is immediately followed by $k+1$ in $v$.
In the former case, it is clear from the definitions of $\L_K(v)$ and $w$ that $k$ is southwest of $k+1$ in the  mountain range description of $w$.
In the latter case, let $(k_j+\ell_j; k_j+\ell_j-1, \ldots, k_j )$ be a consecutive string in $K$ as defined in \eqref{eq_string}.
If $\{k+1, k \}$ is not included in the set $\{k_j+\ell_j, k_j+\ell_j-1, \ldots, k_j \}$ for any $j$, then $k$ is the southwest of $k+1$ in $w$ by an argument similar to the one given above. If $\{k+1, k \}$ is included in $\{k_j+\ell_j, k_j+\ell_j-1, \ldots, k_j \}$ for some $j$, then $k$ is the southwest of $k+1$ in $w$ (see Figures~\ref{picture:definition LKv-1} and \ref{picture:definition LKv-2}).
Thus, we conclude that $w \in \twk$.

(2) For $w \in \twk$ and $u \in \Hop(w)$, we first assert that $\L_K\big(\J_K(u)\big)=u$.

\noindent 
\textbf{Case~(a):} Assume that $\J_K(u) = u$. 
Let $p$ be a peak of $u$. If $u^{-1}(p-1)=u^{-1}(p)+1$, 
then we have $p-1 \not \in K$ by a similar argument to that  in the proof of Lemma~\ref{lemma:kj-1} because  the permutation $w \in \Hop(u)$ belongs to $\twk$.
Hence we have $\L_K(u)=u$ by Definition~\ref{definition:key inverse operator}, which implies $\L_K\big(\J_K(u)\big)=\L_K(u)=u$.

\noindent 
\textbf{Case~(b):} Suppose that $\J_K(u) \neq u$ and consider the effective expression for $\J_K(u)$ in \eqref{eq:effective expression2}.
Here, we recall that $K_j$ is the closed interval $[k_j,k_j+\ell_j-1]$ in $K$ with 
\[
k_1 < k_1+\ell_1 < k_2 < k_2+\ell_2 < \cdots < k_m < k_m+\ell_m
\] 
such that each $\J_{K_j}$ satisfies 
\begin{align*}
\J_{k_j+i}\J_{k_j+i-1} \cdots \J_{k_j}(u_j) \neq \J_{k_j+i-1}\J_{k_j+i-2} \cdots \J_{k_j}(u_j) \ \ \ {\rm for \ all} \ 0 \leq i \leq \ell_j-1
\end{align*}
where $u_1=u$ and $u_j = \J_{K_{j-1}} \cdots \J_{K_1}(u)$ for $2 \leq j \leq m$. 

Let $M$ be the set of strings in $K$ corresponding to $\J_K(u)$ as defined in \eqref{eq_string}. 
Then, it follows from Lemma~\ref{lemma:kj-1} that $(k_j+\ell_j; k_j+\ell_j-1, \ldots, k_j ) \in M$ for all $1 \leq j \leq m$. 
By a similar argument as in Case~(a), we obtain $M=\{(k_j+\ell_j; k_j+\ell_j-1, \ldots, k_j ) \mid 1 \leq j \leq m \}$. 
Putting $v=\J_K(u)$, we get $v_j=u_j$ for all $1 \leq j \leq m$ by Definition~\ref{definition:key inverse operator}. Hence we have $\L_K\big(\J_K(u)\big)=\L_K(v)=v_1=u_1=u$. 

Secondly, we claim that $\J_K(\L_K(v))=v$ for $v \in W(K)$. 
Let \[
\{(k_j+\ell_j; k_j+\ell_j-1, \ldots, k_j ) \mid 1 \leq j \leq m \}
\]
be the set of strings in $K$ corresponding to $v$ as defined in \eqref{eq_string}.
Without loss of generality, we may assume that $k_1 < k_1+\ell_1 < k_2 < k_2+\ell_2 < \cdots < k_m < k_m+\ell_m$.
By Definition~\ref{definition:key inverse operator}, we have a sequence of permutations 
\[
v=v_{m+1}, v_{m}, v_{m-1}, \ldots, v_1=\L_K(v).
\]
It follows from~\eqref{(1)} that the permutation $u\coloneqq v_1(=\L_K(v))$ determines the unique permutation $w \in \twk$ such that $u \in \Hop(w)$.

Let $K_j$ be the closed interval $[k_j,k_j+\ell_j-1]$ for each $1 \leq j \leq m$.
Note that $K_j$ is included in $K$ for all $j$.
We show that $\J_K(u)=\J_{K_m} \cdots \J_{K_1}(u)$.
If $\J_K(u) \neq \J_{K_m} \cdots \J_{K_1}(u)$, then there exists $k \in K$ with $k_{j-1}+\ell_{j-1}< k < k_j$ for some $j$ with $1 \leq j \leq m+1$ such that $\J_{k}(u_{j}) \neq u_{j}$ with $u_j = \J_{K_{j-1}} \cdots \J_{K_1}(u)$. (Here, the inequality $k_{j-1}+\ell_{j-1}< k < k_j$ means $k<k_1$ whenever $j=1$ and $k_m+\ell_m<k$ whenever $j=m+1$, respectively.)
This implies that either Case (i) or Case (ii) below occurs in the mountain range description of $u \in \Hop(w)$. 

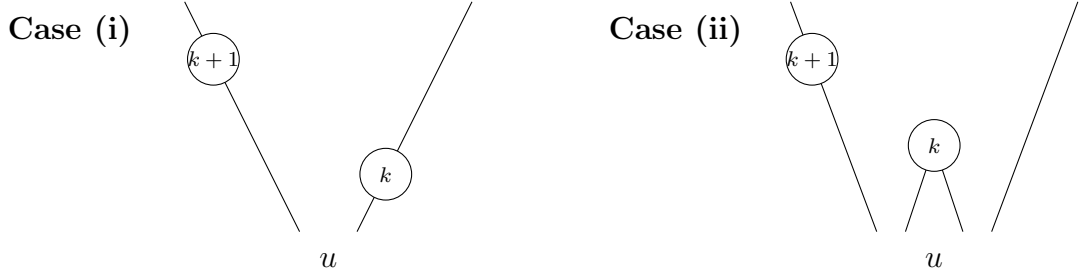
\begin{center}\begin{tikzpicture}[scale=0.38]
\node at (-9,9) {\textbf{Case (i)}};
\draw (-5,10)--(-1,2);
\draw (1,2)--(5,10);
\draw[fill=white] (2,4) circle (0.9cm);
\node at (2,4) {\tiny$k$};

\draw[fill=white] (-4,8) circle (0.9cm);
\node at (-4,8) {\tiny$k+1$};
\node at (0,1) {$u$};

\begin{scope}[xshift=600]
\node at (-9,9) {\textbf{Case (ii)}};
\draw (-5,10)--(-2,2);
\draw (2,2)--(5,10);
\draw (-1, 2)--(0,5)--(1,2);
\draw[fill=white] (0,5) circle (0.9cm);
\draw[fill=white] (-2-9/4,8) circle (0.9cm);
\node at (0,5) {\tiny$k$};
\node at (-2-9/4,8) {\tiny$k+1$};

\node at (0,1) {$u$};
\end{scope}
\end{tikzpicture}
\end{center}

%
%
%
%
%
%
%
%
%
%
%
%
%

By the definition of $u=\L_K(v)$, the positional relationship between $k$ and $k+1$ does not change in the mountain range description of $v$, that is $v^{-1}(k)-v^{-1}(k+1) > 1$. 
However, this is a contradiction since $v \in W(K)$.
Hence, we obtain $\J_K(u)=\J_{K_m} \cdots \J_{K_1}(u)$ which implies that 
\begin{align*}
v&=v_{m+1}=\J_{K_m}(v_m)=\J_{K_m}\J_{K_{m-1}}(v_{m-1})=\cdots=\J_{K_m} \cdots \J_{K_1}(v_1) \\
&=\J_{K_m} \cdots \J_{K_1}(u)=\J_K(u)=\J_K(\L_K(v)).
\end{align*}
Therefore, the result follows. 
\end{proof}

We notice that Proposition \ref{proposition:key} follows immediately from Propositions~\ref{proposition:descent-preserving} and \ref{proposition:bijection}. 

\section{Proof of Theorem~\ref{theorem:main}}
\label{sec_proof_of_main_thm}
Recall from \eqref{eq:2-1} that we have 
\begin{align*} 
\sum_{u \in \Hop(w)} t^{\des(u)} = t^{\pk(w)} (1+t)^{n-1-2\pk(w)} 
\end{align*}
for any $w \in \mathfrak{S}_n$.
By summing the equality above over all $w \in \twk$, we obtain
\begin{align}\label{eq_des_pk_des}
\sum_{w \in \twk} \sum_{u \in \Hop(w)} t^{\des(u)} = \sum_{w \in \twk} t^{\pk(w)} (1+t)^{n-1-2\pk(w)} = \sum_{w \in \twk} t^{\des(w)} (1+t)^{n-1-2\des(w)}. 
\end{align}
Here, the last equality follows from \eqref{eq:hatSn}.
By Propositions~\ref{proposition:h-polynomial PnK} and \ref{proposition:key} the leftmost term of  \eqref{eq_des_pk_des} is equal to  
\begin{align*} 
\sum_{v \in W(K)} t^{\des(v)} = h_{P_n(K)}(t). 
\end{align*}
Therefore, we conclude that  
\begin{align} \label{eq:main} 
h_{P_n(K)}(t) = \sum_{w \in \twk} t^{\des(w)} (1+t)^{n-1-2\des(w)}.
\end{align}

It remains to prove the equality
\begin{align} \label{eq:main2}
|\{w \in \twk \mid \des(w)=j \}| = |\{w \in \hwk \mid \des(w)=j \}|.
\end{align}
We set
$$
K^\ast\coloneqq\{n-k\mid k \in K\},
$$
and for $w \in \S_n$ let $w^\ast$ be the permutation satisfying $w^\ast(i)=n+1-w(n+1-i)$ for all $i \in [n]$.  One can confirm that $\Des(w^\ast)=\{n-i:i \in \Des(w)\}$ and that
$$
\widetilde{W^{K^\ast}}=\{w^\ast\mid w \in \hwk \}.
$$
It follows now from \eqref{eq:main} that
$$
h_{P_n(K^\ast)}(t) = \sum_{w \in \hwk}t^{\des(w)}(1+t)^{n-1-2\des(w)}.
$$
The equality \eqref{eq:main2} will follow once we prove that $P_n(K^\ast)$ and $P_n(K)$ are isomorphic.  To do this, we observe that we may assume that $P_n$ is the convex hull of $\{w(p):w \in S_n\}$, where $p=(p_1,\ldots,p_n)$ satisfies $p_{n+1-i}=-p_i$ for all $i \in [n]$.
Under this assumption, it is straightforward to show that if $w_0$ is the permutation matrix mapping the standard basis vector $e_i$ to $e_{n+1-i}$ for all $i \in [n]$, then $P_n(K)=-w_0(P_n(K^\ast))$.  The claim follows.
This completes the proof of Theorem~\ref{theorem:main}.

\section{Athanasiadis' result} \label{Athanasiadis' result}

In this section we explain Athanasiadis' result, which is an analogue of the result of Foata--Sch\"utzenberger mentioned in the introduction for the graded representation of the symmetric group $\S_n$ on the cohomology of the permutohedral variety. 
We begin with some elementary definitions and known facts from combinatorics.

\subsection{Compositions, partitions, and Young tableaux} \label{sub2}
A {\it composition} of $n$ is a sequence $\alpha=(\alpha_1,\ldots,\alpha_m)$ of positive integers such that $\sum_{i=1}^m\alpha_i=n$.  
For a subset $K \subseteq [n-1]$, the composition $\mu(K)$ of $n$ is defined to be $(|O_1|,\ldots,|O_{n-|K|}|)$ where $O_1,\ldots,O_{n-|K|}$ are $W_K$-orbits in $[n]$
For example, if $n=7$ and $K=\{2,3,4,6\}$ then $O_1=\{1\}$, $O_2=\{2,3,4,5\}$, $O_3=\{6,7\}$, and $\mu(K)=(1,4,2)$. 

The composition $\alpha$ is a {\it partition} of $n$ if $\alpha_i \geq \alpha_{i+1}$ for all $i \in [m-1]$.  Given a finite or infinite sequence $\beta$ of non-negative integers with finitely many nonzero entries, we write $\pr(\beta)$ for the partition obtained by rearranging the entries of $\beta$ so that they appear in nonincreasing order and removing all zeroes.  For example, if $\beta=(4,0,7,5,7,0)$ then $\pr(\beta)=(7,7,5,4)$.
Given a partition $\lambda=(\lambda_1,\ldots,\lambda_\ell)$, the {\it Young diagram} of $\lambda$ consists of $\ell$ rows of squares in the plane, all squares of the same size, each row left-justified against some fixed vertical line, each square in the $(i+1)^{st}$ row having its top edge identified with the bottom edge of a square in row $i$, such that row $i$ contains $\lambda_i$ squares.  For example, the Young diagram of $(3,3,2,1,1)$ appears below. 

\[
\ytableausetup{boxframe=normal, boxsize=1em}
\ydiagram{3,3,2,1,1}
\]

A {\it Young tableau} of shape $\lambda$ is obtained by assigning a positive integer to each square in the Young diagram of $\lambda$.  Such a tableau is {\it semistandard} if the assigned integers increase weakly from left to right in each row and increase strictly from top to bottom in each column.  A semistandard Young tableau of shape $(3,3,2,1,1)$ appears below. 

\[
\ytableausetup{centertableaux}
\begin{ytableau}
2&2&4\\
3&5&5\\
5&6\\
8\\
9
\end{ytableau}
\]

We write $\ssyt_\lambda$ for the set of all semistandard Young tableaux of shape $\lambda$ and $\ssyt_n$ for the set of all semistandard Young tableaux whose underlying Young diagram has $n$ squares.  The {\it content} $c(T)$ of a semistandard Young tableau $T$ is the sequence $(c_1,c_2,\ldots)$ such that $c_i$ is the number of squares in $T$ to which $i$ is assigned.  The tableau above has content $(0,2,1,1,3,1,0,1,1,0,\ldots)$.  Often the content sequence is truncated so that it ends with its last nonzero entry.  We will not distinguish the truncated version from the non-truncated one.

Given an arbitrary semistandard Young tableau $T$, we write $\lambda(T)$ for the partition from which the Young diagram for $T$ is obtained.

Given a partition $\lambda$ of $n$, a {\it standard Young tableau} of shape $\lambda$ is a semistandard Young tableau $T$ of shape $\lambda$ whose content $c(T)$ satisfies $c_i=1$ for $i \in [n]$ (and $c_i=0$ for $i>n$).  For example, the tableau below is a standard Young tableau of shape $(3,2,2)$. 
\[
\ytableausetup{centertableaux}
\begin{ytableau}
1&3&7\\
2&5\\
4&6
\end{ytableau}
\]
We write $\syt_\lambda$ for the set of all standard Young tableaux of shape $\lambda$ and $\syt_n$ for the set of all standard Young tableaux whose underlying Young diagram has $n$ squares.

Given $T \in \syt_n$, a {\it descent} of $T$ is any $i \in [n-1]$ such that $i$ appears in a higher row of $T$ than does $i+1$.  We write $\Des(T)$ for the set of descents of $T$ and $\des(T)$ for $|\Des(T)|$.  For example, with $T$ as above, $\Des(T)=\{1,3,5\}$ and $\des(T)=3$.  Double, initial and final descents are defined as they were for words.  We write $\widetilde{\syt_\lambda}$ and $\widetilde{\syt_n}$ for the respective subsets of $\syt_\lambda$ and $\syt_n$ consisting or tableaux with no double descent and no initial descent. Similarly, we denote by $\widehat{\syt_\lambda}$ and $\widehat{\syt_n}$ the respective subsets of $\syt_\lambda$ and $\syt_n$ consisting or tableaux with no double descent and no final descent.

\subsection{The Robinson--Schensted--Knuth (RSK) algorithm} 
We will not describe the RSK algorithm here.  Rather, in addition to providing an example below, we will describe its key consequences.  The interested reader can find details in \cite[Chapter 7]{StanleyEC}.  With the notation used in (7.38) on page 318 of \cite{StanleyEC}, we assume at all times that $i_k=k$ for all $k \in [m]$.  Under this assumption the algorithm produces from a word $v$ of length $n$ a pair $(P_v,Q_v) \in \ssyt_n \times \syt_n$ satisfying $\lambda(P_v)=\lambda(Q_v)$ along with two other key properties.

\begin{proposition} \label{RSKproperties}
For every positive integer $n$, there is a bijection 
$$
\rsk:\Z_{>0}^n \rightarrow \{(P,Q) \in \ssyt_n \times \syt_n \mid \lambda(P)=\lambda(Q)\}
$$
such that given $\rsk(v)=(P_v,Q_v)$, both
\begin{itemize}
\item[(A)] if $P_v$ has content $c$, then $c_i=m_i(v)$ for all $i$, and
\item[(B)] $\Des(Q_v)=\Des(v)$
\end{itemize}
hold.
\end{proposition}

For example, applying RSK to the word $v=23132$, we obtain the following sequence of pairs of tableaux, terminating with $(P_v,Q_v)$.  

\[
\left(\begin{ytableau} 2 \end{ytableau},
\begin{ytableau} 1 \end{ytableau} \right) \to
\left(\begin{ytableau} 2&3 \end{ytableau},
 \begin{ytableau} 1&2 \end{ytableau} \right) \to
\left(\begin{ytableau} 1&3\\2 \end{ytableau},
 \begin{ytableau} 1&2\\3 \end{ytableau} \right) \to
\left(\begin{ytableau} 1&3&3\\2 \end{ytableau},
 \begin{ytableau} 1&2&4\\3 \end{ytableau} \right) \to
\left(\begin{ytableau} 1&2&3\\2&3 \end{ytableau},
 \begin{ytableau} 1&2&4\\3&5 \end{ytableau} \right). 
\]

We observe that $Q_v \in \widetilde{\syt_n}$ (resp. $Q_v \in \widehat{\syt_n}$) if and only if $v \in \widetilde{\Z_{>0}^n}$ (resp. $v \in \widehat{\Z_{>0}^n}$) and that $P_v \in \syt_n$ if and only if $v \in \S_n$.  The fact that (A) holds follows immediately from the definition of the RSK algorithm.  The fact that (B) holds follows from the discussion following (7.43) on page 321 along with Lemma 7.23.1 in \cite{StanleyEC}.  

\subsection{Modules for $\C[\S_n]$, the Frobenius characteristic, and Kostka numbers}
The irreducible modules of $\C[\S_n]$ are parameterized by partitions of $n$.  We write $S^\lambda$ for the module associated to the partition $\lambda$.  Every $\C[\S_n]$-module is a direct sum of simple modules, and each $S^\lambda$ is a summand in some permutation module $M^\mu$.  Here $\mu=\mu(K)$ for some $K \subseteq [n-1]$ and $M^\mu$ has a $\C$-basis indexed by the cosets of $W_K$ in $\S_n$.  The elements of $\S_n$ permute this basis as they permute the cosets by translation.  The {\it Kostka number} $\kostka_{\lambda,\mu}$ is the number of semistandard Young tableaux of shape $\lambda$ and content $\mu$.  The first equality in the next theorem is well-known (see for example \cite[Theorem 7.18.7]{StanleyEC} and the second follows directly from Frobenius reciprocity (see for example \cite[Lemma~5.2]{Isaacs}).

\begin{theorem} \label{kostka}
The Kostka number $\kostka_{\lambda,\mu}$ is equal to the number of summands isomorphic to $S^\lambda$ in a direct sum decomposition of $M^\mu$ into simple submodules, which is equal in turn to the dimension $\dim (S^\lambda)^{W_K}$ of the subspace of $S^\lambda$ consisting of vectors fixed by a parabolic subgroup $W_K$ satisfying $\mu(K)=\mu$.  
\end{theorem}

The {\it Frobenius characteristic} $\ch$ is a map from the direct sum of the spaces of virtual modules for the group algebras $\C[\S_n]$ ($n \in \Z_{\geq 0}$) to the algebra of symmetric functions.  We do not need to know anything about this map other than that for each $n$ it determines an isomorphism between the space of virtual $\C[\S_n]$-modules and the space of homogeneous degree $n$ symmetric functions.  All of this is described (using character-theoretic language) in \cite[Section 7.18]{StanleyEC}.

\subsection{Graded $\S_n$-representation formula for gamma vectors} 
Consider the permutohedral variety $X_n$, which is the toric variety associated with the normal fan $\Delta_n$ of a permutohedron $P_n$.
Since the fan $\Delta_n$ is invariant under the reflection action of $\S_n$, this induces an action of $\S_n$ on the permutohedral variety $X_n$, which determines a representation of $\S_n$ on $H^\ast(X_n;\C)$.
This representation is palindromic as noted in \cite[Corollary 4.4]{Stanley1} and \cite[Proposition 12]{Stembridge1}, which implies that there exist virtual modules $R_{n,j}$ ($0 \leq j \leq \lfloor \frac{n-1}{2} \rfloor$) such that
\begin{equation} \label{rdef}
\sum_{k=0}^{n-1}H^{2k}(X_n;\C)t^k \cong_{\S_n} \sum_{j=0}^{\lfloor \frac{n-1}{2} \rfloor}R_{n,j}t^j(1+t)^{n-1-2j},
\end{equation}
where $\cong_{\S_n}$ denotes isomorphism of $\C[\S_n]$-modules. 
Athanasiadis proved the following result.

\begin{theorem}$($\cite[Corollary 2.41]{Athanasiadis}$)$ \label{theorem:Athanasiadis}
For $0 \leq j \leq \lfloor \frac{n-1}{2} \rfloor$,
$$
R_{n,j} \cong_{\S_n} \sum_{\substack{Q \in \widehat{\syt_n} \\ {\des(Q)=j}}}S^{\lambda(Q)}= \sum_{\substack{Q \in \widetilde{\syt_n} \\ {\des(Q)=j}}}S^{\lambda(Q)}.
$$
Equivalently,
\begin{align*}
\sum_{k=0}^{n-1}H^{2k}(X_n;\C)t^k &\cong_{\S_n} \sum_{Q \in \widehat{\syt_n}}S^{\lambda(Q)}t^{\des(Q)}(1+t)^{n-1-2\des(Q)}\\
&= \sum_{Q \in \widetilde{\syt_n}}S^{\lambda(Q)}t^{\des(Q)}(1+t)^{n-1-2\des(Q)}.
\end{align*}
\end{theorem}

For the reader's convenience, we here give a proof of Theorem~\ref{theorem:Athanasiadis}.
We begin with the result of Gessel (cf. \cite[Theorem~7.3]{SWEQF}).  Given $v \in \widehat{\Z_{>0}^n}$, we define $\x^v\coloneqq\prod_{i \geq 1}x_i^{m_i(v)}$.  So, the monomial $\x^v$ tells us how many times each letter $i$ appears in the word $v$.  As stated in \cite{SWEQF},
\begin{equation} \label{gessel}
\sum_{k=0}^{n-1}\ch(H^{2k}(X_n;\C))t^k=\sum_{v \in \widehat{\Z_{>0}^n}}\x^vt^{\des(v)}(1+t)^{n-1-2\des(v)}.
\end{equation}
It follows from \eqref{gessel} and the definition of $R_{n,j}$ in \eqref{rdef} that for $j \in [n-1]$
\begin{equation} \label{gessel2}
\ch(R_{n,j})=\sum_{\substack{v \in \widehat{\Z_{>0}^n} \\ \des(v)=j}}\x^v.
\end{equation}

Given a Young tableau $T$, we write $\x^T$ for the monomial $\prod_{i \geq 1}x_i^{m_i(T)}$, where $m_i(T)$ is the number of boxes in $T$ filled with $i$.   For any partition $\lambda$ of $n$,
$$
\ch(S^\lambda)=\sum_{T \in \ssyt_\lambda}\x^T,
$$
see for example \cite[Definition 7.10.1 and (7.86)]{StanleyEC}.  Therefore, for $Q \in \syt_n$,
\begin{equation} \label{schur}
\ch(S^{\lambda(Q)})=\sum_{P \in \ssyt_{\lambda(Q)}}\x^P.
\end{equation}

We conclude now that for $j \in [n-1]$
\begin{eqnarray*}
\sum_{\substack{Q \in \widehat{\syt_n} \\ \des(Q)=j}}\ch(S^{\lambda(Q)}) =  \sum_{\substack{(P,Q) \in \ssyt_n \times \widehat{\syt_n} \\ \lambda(P)=\lambda(Q) \\ \des(Q)=j}}\x^P = \sum_{\substack{v \in \widehat{\Z_{>0}^n} \\ \des(v)=j}} \x^v  = \ch(R_{n,j}),
\end{eqnarray*}
the first equality following from \eqref{schur}, the second from Proposition \ref{RSKproperties} and the third from~\eqref{gessel2}.  As $\ch$ is an isomorphism, the first half of Theorem \ref{theorem:Athanasiadis} follows.
It remains to show the following equality 
\begin{align} \label{eq:hattildetableaux}
\sum_{\substack{Q \in \widehat{\syt_n} \\ {\des(Q)=j}}}S^{\lambda(Q)} = \sum_{\substack{Q \in \widetilde{\syt_n} \\ {\des(Q)=j}}}S^{\lambda(Q)}.
\end{align}
For this, we consider the evacuation $Q^\ast$ for $Q \in \syt_n$. 
We will not describe the evacuation, but we will describe its key consequences. 
We refer the reader to \cite[Appendix~A.1]{Fulton} for the details. 
The evacuation $Q^\ast$ for $Q \in \syt_n$ is the standard tableau of entries $i^\ast \coloneqq n+1-i \ (i \in [n])$ with the same shape $\lambda(Q)$. 
For $w \in \S_n$, let $w^\ast$ be the permutation defined by $w^\ast(i)=n+1-w(n+1-i)$ for all $i \in [n]$. 
It is known that if $\rsk(w)=(P,Q)$, then $\rsk(w^\ast)=(P^\ast,Q^\ast)$.
This property implies that 
\begin{align*}
\Des(Q^\ast) = \Des(w^\ast) = \{ (i+1)^\ast \mid i \in \Des(w) \} = \{ (i+1)^\ast \mid i \in \Des(Q) \}.
\end{align*}
The equality above shows that the assigement by sending $Q$ to $Q^*$ gives a one-to-one correspondence
\begin{align*}
\{Q \in \widehat{\syt_n} \mid \des(Q)=j \} \rightarrow \{ Q \in \widetilde{\syt_n} \mid \des(Q) = j \},
\end{align*}
which yields \eqref{eq:hattildetableaux}.

\section{A relation between our theorem and Athanasiadis' result} \label{relation}
The aim of this section is to explain the relation between our main theorem and Athanasiadis' result. 
It is not surprising that Theorem \ref{theorem:main} follows directly from Theorem \ref{theorem:Athanasiadis} and vice versa.  
We begin with a result from \cite{HMSS}.
For a subset $K \subseteq [n-1]$, we write $X_n(K)$ for the toric orbifold associated with the partitioned permutohedron $P_n(K)$. 

\begin{theorem}[\cite{HMSS}, Theorem 1.1] \label{hmss}
For any subset $K \subseteq [n-1]$, the cohomology ring $H^\ast(X_n(K);\C)$ is isomorphic with the ring $H^\ast(X_n;\C)^{W_K}$ of $W_K$-fixed elements of $H^\ast(X_n;\C)$.
In particular, we obtain
\begin{eqnarray} \label{eq:hmss}
h_{P_n(K)}(t) = \sum_{\ell=0}^{n-1} \dim H^{2\ell}(X_n(K);\C)t^\ell  = \sum_{\ell=0}^{n-1} \dim(H^{2\ell}(X_n;\C)^{W_K})t^\ell. 
\end{eqnarray}
\end{theorem}

Assume that $K \subseteq [n-1]$ with $\mu(K)=(\mu_1,\ldots,\mu_m)$, a composition of $n$. We define
$$
\widehat{[m]^n}(K)\coloneqq\{v \in \widehat{\Z_{>0}^n} \mid m_i(v)=\mu_i \mbox{ for all } i \in [m]\}.
$$
So, the set $\widehat{[m]^n}(K)$ consists of all words $v$ of length $n$ with entries in $[m]$ such that $v$ has no double descent, no final descent, and $\mu_i$ letters equal to $i$ for each $i \in [m]$ (and no letter greater than $m$).  We observe that $\rsk$ maps $\widehat{[m]^n}(K)$ onto the set of pairs $(P,Q)$ from $\ssyt_n \times \widehat{\syt_n}$ such that $\lambda(P)=\lambda(Q)$ and $P$ has content $\mu(K)$.  Therefore,
\begin{equation} \label{hatmn}
\sum_{v \in \widehat{[m]^n}(K)}t^{\des(v)}(1+t)^{n-1-2\des(v)}=\sum_{Q \in \widehat{\syt_n}}\kostka_{\lambda(Q),\mu(K)}t^{\des(Q)}(1+t)^{n-1-2\des(Q)}.
\end{equation}

Here, one can see that there is a descent-preserving bijection from $\widehat{[m]^n}(K)$ to $\widehat{W^K}$.  Given $v \in \widehat{[m]^n}(K)$ and $i \in [m]$, define
$$
X_i(v)\coloneqq\{j \in [n]\mid v(j)=i\}
$$
and let $\phi(v) \in \S_n$ be the permutation such that if $j$ is the $k^{th}$ smallest element of $X_i(v)$ then 
$$
\phi(v)_j=k+\sum_{\ell=1}^{i-1}\mu_\ell.
$$
For example, if $n=7$ and $K=\{1,2,4,6\}$, we have $\mu(K)=(3,2,2)$. For an element $v=3231211 \in \widehat{[m]^n}(K)$, it follows that $\phi(v)=6471523 \in \widehat{W^K}$. 

It is straightforward to confirm that $\phi$ is the desired bijection.
Hence, the left hand side of \eqref{hatmn} is described as $\sum_{w \in \hwk}t^{\des(w)}(1+t)^{n-1-2\des(w)}$.
On the other hand, the right hand side of \eqref{hatmn} is $\sum_{Q \in \widehat{\syt_n}}\dim(S^{\lambda(Q)})^{W_K}t^{\des(Q)}(1+t)^{n-1-\des(Q)}$ by Theorem \ref{kostka}.
Thus, we can rewrite \eqref{hatmn} as 
\begin{equation} \label{hatmn2}
\sum_{w \in \hwk}t^{\des(w)}(1+t)^{n-1-2\des(w)} = \sum_{Q \in \widehat{\syt_n}}\dim(S^{\lambda(Q)})^{W_K}t^{\des(Q)}(1+t)^{n-1-\des(Q)}.
\end{equation}
Therefore, we conclude from \eqref{eq:hmss} and \eqref{hatmn2} that the equation 
\[
h_{P_n(K)}(t) = \sum_{w \in \hwk}t^{\des(w)}(1+t)^{n-1-2\des(w)}
\] 
is equivalent to
\[
\sum_{\ell=0}^{n-1} \dim(H^{2\ell}(X_n;\C)^{W_K})t^\ell = \sum_{Q \in \widehat{\syt_n}}\dim(S^{\lambda(Q)})^{W_K}t^{\des(Q)}(1+t)^{n-1-\des(Q)}.
\]

If we pick for each partition $\nu$ of $n$ some $K_\nu$ such that $\pr(\mu(K_\nu))=\nu$, then any $\C[S_n]$-module $M$ is determined up to isomorphism by the dimensions of all of the spaces of $W_{K_\nu}$-invariant vectors.  (This is well-known and follows from, for example \cite[Propositions 7.10.5 and 7.18.7]{StanleyEC}.)
Now we know that Theorem \ref{theorem:main} follows quickly from Theorem \ref{theorem:Athanasiadis} and vice versa.

\end{document}